\documentclass[12pt,a4paper]{article}

\usepackage{latexsym}
\usepackage[pdftex]{graphicx}
\usepackage{amssymb}
\usepackage{amsmath}
\usepackage{float} 

\usepackage[pdftex]{graphicx}
\newtheorem{theorem}{Theorem}[section]

\newtheorem{remark}{Remark}[section]

\newtheorem{conjecture}{Conjecture}[section]

\begin{document}
\begin{center}
{\bf \LARGE On the composition and decomposition of positive linear operators III: A non-trivial decomposition of the Bernstein operator}\\[0.25cm]
{\bf \LARGE -- Technical report --}
\end{center}
{\large \bf Heiner Gonska, Margareta Heilmann,  \framebox{Alexandru Lupa\c{s}}, Ioan Ra\c{s}a}
\vspace{0.5cm}
\begin{abstract}
The central problem in this technical report is the question if the classical Bernstein operator can be decomposed into 
nontrivial building blocks where one of the factors is the genuine Beta operator introduced by 
 M\"uhlbach \cite{Mu1970} and  Lupa\c{s} \cite{Lu1972}. 
\\
We collect several properties of the Beta operator such as injectivity, the eigenstructure and the images of the monomials
under its inverse.
Moreover, we give a decomposition of the form $B_n = \bar{\mathbb{B}}_n  \circ F_n $ where $F_n$ is a nonpositive linear operator having quite interesting properties. 
We study the images of the monomials under $F_n$, its moments and various representations.
Also an asymptotic formula of Voronovskaya type for polynomials is given and a connection with a conjecture
of Cooper and Waldron \cite{CoWa2000} is established.
In an appendix numerous examples illustrate the approximation behaviour of $F_n$ in comparison to $B_n$.
\end{abstract}
\\[0.25cm]
\section{Introduction}
\label{intro}
Here we continue our previous research on the composition of positive linear operators and on linear operators in general, thus emphasizing the fact again that a functional-analytic point of view onto the problem is useful. Our report is a continuation of \cite {Go2000} and \cite{GoRa2010} where related problems were considered.

The present report is motivated by a discussion between the late Alexandru Lupa\c{s} and the first author 
which took place in Sibiu in late December 2006.
The topic of this discussion was the question if the classical Bernstein operator 
$$
	B_n(f;x) = \sum_{k=0}^n f \left ( \frac{k}{n} \right ) {n \choose k} x^k (1-x)^{n-k} ,
$$
$f \in C[0,1]$, $x \in [0,1]$, can be decomposed into simpler positive building blocks.
More precisely, the problem was if there are non-trivial positive linear operators
$P$ and $Q$ such that $B_n = P \circ Q$.
We had some preliminary results then, and it was intended to eventually publish a joint paper dealing with this topic. This is mentioned in the obituary \cite{GaGoKaPaStVe2009} indicating one reason for the long delay in further investigating the problem.

One of our candidates for the factors $P$ and $Q$ were certain Beta-type operators introduced by M\"uhlbach in \cite{Mu1970} and further investigated by him in \cite{Mu1972} and by Lupa\c{s} in \cite{Lu1972}. 
These mappings are given for $f \in C[0,1]$, $x \in [0,1]$ by
$$
	\bar{\mathbb{B}}_n (f;x) =
	\begin{cases}
	f(0) &,\ x=0 ,\\
	\displaystyle \frac{1}{B(nx,n(1-x))} \int_0^1 t^{nx-1} (1-t)^{n(1-x)-1} f(t) dt &, \ x \in (0,1) , \\
	f(1) &, \ x = 1 .
	\end{cases}
$$
Here $B( \cdot , * )$ is the Beta function. The $\bar{\mathbb{B}}_n$ are positive linear endomorphisms of $C[0,1]$;
they reproduce linear functions and have second moments smaller than those of the Bernstein operators. More
precisely, see \cite[Satz 2.28]{Lu1972},
$$
	\bar{\mathbb{B}}_n ((e_1 -x)^2 ; x) = \frac{x(1-x)}{n+1} \leq \frac{x(1-x)}{n} 
	= B_n ((e_1 -x)^2 ; x) .
$$
Moreover, it is known from \cite{AdBaCa1993} and \cite{AtRa2007} that $\bar{\mathbb{B}}_n$ preserves monotonicity
and (ordinary) convexity.

It is known that if one composes two positive linear operators $P$ and $Q$, both reproducing linear functions, then for the second moment of the product operator one has
$$
	( P \circ Q )  ((e_1 -x)^2 ; x) = P^u (Q  ((e_1 -u)^2 ; u); x) + P  ((e_1 -x)^2 ; x) .
$$
Here the superscript in $P^u$ indicates that the operator $P$ is applied to functions in the variable $u$.

Putting $P = \bar{\mathbb{B}}_n$ the question then was if there is another positive linear operator $Q$ such that 
$ \bar{\mathbb{B}}_n \circ Q = B_n$ and in particular,
\begin{eqnarray*}
	(  \bar{\mathbb{B}}_n \circ Q )   ((e_1 -x)^2 ; x)
	& = &
	B_n   ((e_1 -x)^2 ; x)
	\\
	& = &
	\frac{x(1-x)}{n}
	\\
	& = &
	  \bar{\mathbb{B}}_n^u (Q  ((e_1 -u)^2 ; u); x) +  \bar{\mathbb{B}}_n  ((e_1 -x)^2 ; x)
	\\
	& = &
	 \bar{\mathbb{B}}_n^u (Q  ((e_1 -u)^2 ; u); x) + \frac{x(1-x)}{n+1} .
\end{eqnarray*}
Natural candidates for $Q$ are operators of the form
$$
	Q (f;x) = \sum_{k=0}^n f \left ( \frac{k}{n} \right ) r_{n,k} (x) ,
$$
with $r_{n,k} \geq 0$, $x \in [0,1]$, $0 \leq k \leq n$, so that
$$
	(  \bar{\mathbb{B}}_n \circ Q )  (f;x) =  \sum_{k=0}^n f \left ( \frac{k}{n} \right )  \bar{\mathbb{B}}_n (r_{n,k},x) 
$$
would become the Bernstein operator if $r_{n,k} $ could be chosen in a way such that
$$
	 \bar{\mathbb{B}}_n (r_{n,k},x) = b_{n,k} (x) :=  {n \choose k} x^k (1-x)^{n-k} , \ x \in [0,1] , \ 0 \leq k \leq n.	
$$
The first (unpublished) attempt made used piecewise linear interpolation 
$$
	 S_{\Delta_n} \, : \, C[0,1] \longrightarrow C[0,1] \text{ at }
	0 , \frac{1}{n} , \dots ,  \frac{k}{n}, \dots , \frac{n-1}{n} , 1
$$
which can explicitely be described as 
$$
	 S_{\Delta_n} (f;x) = \frac{1}{n} \sum_{k=0}^n 
	\left [ \frac{k-1}{n} ,  \frac{k}{n} ,  \frac{k+1}{n} ; |\alpha -x| \right ]_{\alpha} 
	 f \left ( \frac{k}{n} \right ) ,
$$
where $[a,b,c ;f] = [a,b,c; f(\alpha) ]_{\alpha}$ denotes the divided difference of a function
$f \, : \, D \longrightarrow  \mathbb{R}$ on the (distinct knots) $\{a,b,c\} \subset D$, with respect to $\alpha$.
$ S_{\Delta_n}$ is also a positive linear operator reproducing linear functions and preserving monotonicity and 
convexity/concavity. Moreover, it is of the appropriate form and hence it made sense to consider 
${\mathbb{G}}_n:=  \bar{\mathbb{B}}_n \circ S_{\Delta_n} $, that is,
$$
	{\mathbb{G}}_n \, : \, C[0,1] \longrightarrow C[0,1] ,
$$
where
$$
	{\mathbb{G}}_n (f;0) =  S_{\Delta_n} (f;0) = f(0)  , \,
	{\mathbb{G}}_n (f;1)= S_{\Delta_n} (f;1) = f(1),
$$
and for $x \in (0,1)$
$$
	{\mathbb{G}}_n (f;x)=  \frac{1}{B(nx,n(1-x))} \int_0^1 t^{nx-1} (1-t)^{n(1-x)-1} S_{\Delta_n} (f;t) dt .
$$
${\mathbb{G}}_n $ is again positive and linear. As the composition of two operators preserving monotonicity and convexity, ${\mathbb{G}}_n $ also has these properties.

For a convex function $g$ it is well-known that $g \leq B_ng$.
Now if $f \in C[0,1]$ is  convex, then this is also true for $ S_{\Delta_n} f$, so that
$$
	 f \leq  S_{\Delta_n} f \leq B_n ( S_{\Delta_n} f ) = B_n f ,
$$
implying
$$
	 \bar{\mathbb{B}}_n f \leq ( \bar{\mathbb{B}}_n \circ  S_{\Delta_n}) f
	={\mathbb{G}}_n f  \leq(  \bar{\mathbb{B}}_ n \circ  B_n )  f = L_n f,
$$
where $L_n$ is a special case of the Stancu operator introduced in
\cite {St1968}, namely for the case $\alpha = \frac{1}{n}$.
In particular,
\begin{eqnarray*}
	\bar{\mathbb{B}}_n    ((e_1 -x)^2 ; x)
	& = &
	\frac{x(1-x)}{n+1}
	\\
	& \leq &
	{\mathbb{G}}_n ((e_1 -x)^2 ; x) 
	\\
	& \leq &
	L_n  ((e_1 -x)^2 ; x)
	\\
	& = &
	  \frac{2x(1-x)}{n+1} .
\end{eqnarray*}
More generally, for $j \in \mathbb{N}_0,$
$$
	\bar{\mathbb{B}}_n    ((e_1 -x)^{2j} ; x)
	\leq {\mathbb{G}}_n ((e_1 -x)^{2j} ; x) 
	\leq L_n  ((e_1 -x)^{2j} ; x) .
$$
The latter inequalities can be used to give estimates for the degree of approximation by ${\mathbb{G}}_n$, but
we will not further discuss this here. Since the second moments of both ${\mathbb{G}}_n$ and $B_n$ lie
between $\frac{x(1-x)}{n+1}$ and $\frac{2x(1-x)}{n+1}$, there still is a chance that 
${\mathbb{G}}_n=B_n$. However, in the next section we will show that  ${\mathbb{G}}_2\not= B_2$.
Moreover, in Section \ref{operators} it will be proved that there is {\bf no} positive linear operator
$Q \, : \, C[0,1] \longrightarrow \Pi_n$ such that $B_n =\bar{\mathbb{B}}_n \circ Q$.
We will also show that it is impossible to write $B_n = L \circ  S_{\Delta_n}$ for a large class of positive integral operators.

But these negative results do not exclude the possibility that there are non-trivial decompositions $B_n = P \circ Q$
with $P \not= \bar{\mathbb{B}}_n $ or $Q \not= S_{\Delta_n}$.
But if one insists in the choice  $P = \bar{\mathbb{B}}_n $, then we are necessarily led to certain non-positive operators $F_n$ which will be mainly investigated in this report starting from Section \ref{operators}.

\section{Two negative results}
\label{negative}
We shall prove that ${\mathbb{G}}_2\not= B_2$. Indeed,
$$
	{\mathbb{G}}_2 f = \sum_{i=0}^2 f \left ( \frac{i}{2} \right )   \bar{\mathbb{B}}_2 u_i , \, f \in C[0,1 ,]
$$
where $u_i \in C[0,1]$ is the piecewise linear function with $u_i  \left ( \frac{j}{2} \right ) = \delta_{ij}$,
$i,j \in \{0,1,2\}$.

Suppose that ${\mathbb{G}}_2= B_2$. Then $ \bar{\mathbb{B}}_2 u_i = b_{2,i}$, $i=0,1,2$.
In particular, $ \bar{\mathbb{B}}_2 u_2 (x) = x^2$, $x \in [0,1]$, which leads to
$$
	\frac{\int_{\frac{1}{2}}^1 t^{2x-1}(1-t)^{1-2x}(2t-1) dt}{B(2x,2(1-x))} 
	= x^2 , \, x \in (0,1) .
$$
For $x = \frac{1}{4}$ we get
\begin{equation}
\label{eq-2-1}
	\frac{   \int_{\frac{1}{2}}^1 t^{-\frac{1}{2}}(1-t)^{\frac{1}{2}}(2t-1) dt   }{\int_{0}^1 t^{-\frac{1}{2}}(1-t)^{\frac{1}{2}} dt} 
	= \frac{1}{16} .
\end{equation}
On $(0,1)$,
$$
	 \int t^{-\frac{1}{2}}(1-t)^{\frac{1}{2}}(2t-1) dt = 
	-\frac{1}{4} \left \{ (6-4t) \sqrt{t(1-t)}+ \arcsin{(2t-1)} \right \}
$$
and
$$
	\int t^{-\frac{1}{2}}(1-t)^{\frac{1}{2}} dt =
	 \sqrt{t(1-t)} + \frac{1}{2}  \arcsin{(2t-1)}
$$
Now (\ref{eq-2-1}) becomes
$$
	\frac{\frac{1}{2}-\frac{\pi }{8}}{\frac{\pi }{2}}= \frac{1}{16},
$$
i.e. $\pi = \frac{16}{5}$, a contradiction.
This proves  ${\mathbb{G}}_2\not= B_2$. 

The next considerations show that it is not possible to write  $B_n = L \circ  S_{\Delta_n}$ for a large class
 of integral operators.
The operator $S_{\Delta_n} \, : \, C[0,1] \longrightarrow C[0,1] $ can be described as in Section \ref{intro}, or as
$$
	S_{\Delta_n}  f (x) = \sum_{i=0}^n  f \left ( \frac{i}{n} \right ) u_{n,i} (x) , \,
	f \in C[0,1], \, x \in [0,1],
$$
where $  u_{n,i} \in  C[0,1]$ are piecewise linear functions such that 
$u_{n,i}  \left ( \frac{j}{n} \right ) = \delta_{ij}$, $i,j = 0, \dots , n$.

Let  $L \, : \, C[0,1] \longrightarrow C[0,1]$ be an integral operator,
$$
	L (f;x) := \int_0^1 K(x,t) f(t) dt , \, f \in C[0,1], \, x \in [0,1],
$$
where the kernel $K$ is non-negative on $[0,1]^2$ and $K(x, \cdot ) \in L_1 [0,1]$ for all $x \in [0,1]$. 
\\
We shall prove that
$L \circ S_{\Delta_n} \not= B_n$, $n \geq 2$.

Suppose that for a given $n \geq 2$ we have $L \circ S_{\Delta_n} = B_n$. Then
$$
	 \sum_{i=0}^n L( u_{n,i}; x)  f \left ( \frac{i}{n} \right ) =
	 \sum_{i=0}^n b_{n,i} (x)  f \left ( \frac{i}{n} \right )  , \, f \in C[0,1],
$$
which entails
$$
	 L( u_{n,i}; x) =  b_{n,i} (x) , \, x \in [0,1] , \, i = 1, \dots , n.
$$
In particular,  $  L( u_{n,i}; 0) = 0 $, $  i = 1, \dots , n $, and so we get
$$
	 \int_0^1 K(0,t)  u_{n,i} (t) dt =  L( u_{n,i}; 0)=0 , \, i = 1 , \dots , n .
$$
It follows that 
$$
	 \int_0^1 K(0,t) \left \{ \sum_{i=1}^n u_{n,i} (t) \right \} dt = 0.
$$
But $ \sum_{i=1}^n u_{n,i} (t) = 1 - u_{n,0} (t) > 0$, for all $t \in (0,1]$.
We deduce that $K(0, \cdot ) = 0$ a.e. on $[0,1]$,
and so 
\begin{equation}
\label{eq-2-2}
L(e_0 ; 0) = \int_0^1 K(0,t) dt = 0 .
\end{equation}
On the other hand,
$$
	L(e_0 ; 0) = ( L \circ   S_{\Delta_n}) (e_0 ; 0) = B_n (e_0 ; 0) = 1 ,
$$
which contradicts (\ref{eq-2-2}). Thus, in fact,  $ L \circ  S_{\Delta_n} \not= B_n$. 

\section{Injectivity of $\bar{\mathbb{B}}_n$}
\label{injectiv}
In this section we will prove that $\bar{\mathbb{B}}_n \, : \, C[0,1] \longrightarrow C[0,1] $ is injective. This fact has the consequence that the operators $F_n$ considered below are the only endomorphisms on $C[0,1]$ allowing the
decomposition
$$
	B_n = \bar{\mathbb{B}}_n \circ F_n ,
$$
meaning that any other endomorphism $Q$ with $B_n = \bar{\mathbb{B}}_n \circ Q$ is necessarily equal to $F_n$.
\begin{theorem}
\label{thm1}
$\bar{\mathbb{B}}_n   \, : \, C[0,1] \longrightarrow C[0,1] $ is injective, $n \geq 1$.
\end{theorem}
{\bf Proof.}
Let $f \in C[0,1]$, $\bar{\mathbb{B}}_n (f,x) =0$, $x \in [0,1]$. Then
$$
	\int_0^1 t^{nx-1} (1-t)^{n(1-x)-1} f(t) dt = 0 , \, \forall x \in (0,1) .
$$
Setting $t:= \frac{{\rm{e}}^u}{1+{\rm{e}}^u}$ and $g(u) := f \left (  \frac{{\rm{e}}^u}{1+{\rm{e}}^u} \right )$, $u \in \mathbb{R}$, we get
\begin{equation}
\label{eq-3-1}
	\int_{\mathbb{R}} {\rm{e}}^{nux} \frac{1}{(1+{\rm{e}}^u)^n} g(u) du = 0, \, x \in (0,1) .
\end{equation}
Obviously $g \in C_b (\mathbb{R})$, i.~e., $\sup_{x \in \mathbb{R}} |g(x)| = \| f \|_{\infty} <  \infty$.
\\
Let $0 < l < \frac{1}{2}$ be fixed, and let $x \in \left [ \frac{1}{2}-l, \frac{1}{2} +l \right ]$.
\\
a) For $u \in [0,\infty)$ and $k \in \mathbb{N}_0$ we have
\begin{eqnarray*}
	{\rm e}^{nux} \frac{1}{(1 + {\rm e}^u)^n} (nu)^k |g(u)|
	& \leq &
	n^k \|f\|_{\infty} u^k \frac{{\rm e}^{nu(\frac{1}{2}+l)}}{{\rm e}^{nu}}
	\\
	& = &
	n^k \|f\|_{\infty} u^k {\rm e}^{(l-\frac{1}{2})nu} .
\end{eqnarray*}
From \cite[p. 708]{Fi1966} it follows that for all $k \in \mathbb{N}_0$ the integral 
$$
	\int_{0}^{\infty} {\rm{e}}^{nux} (nu)^k \frac{1}{(1+{\rm{e}}^u)^n} g(u) du 
$$
is convergent, uniformly with respect to $x \in \left [ \frac{1}{2}-l, \frac{1}{2} +l \right ]$.
\\
b) If $u \in ( -\infty , 0]$ and $k \in \mathbb{N}_0$, then
$$
	{\rm e}^{nux} \frac{1}{(1 + {\rm e}^u)^n} | (nu)^k g(u)|
	\leq 
	n^k \|f\|_{\infty}(- u)^k {\rm e}^{(\frac{1}{2}-l)nu} .
$$
Thus the integral 
$$
	\int_{-\infty }^{0} {\rm{e}}^{nux} (nu)^k \frac{1}{(1+{\rm{e}}^u)^n} g(u) du 
$$
is convergent, uniformly with respect to $x \in \left [ \frac{1}{2}-l, \frac{1}{2} +l \right ]$.
\\
From a) and b) we conclude that for all $k \in \mathbb{N}_0$ the integral
$$
	\int_{-\infty }^{\infty} {\rm{e}}^{nux} (nu)^k \frac{1}{(1+{\rm{e}}^u)^n} g(u) du 
$$
is convergent, uniformly with respect to $x \in \left [ \frac{1}{2}-l, \frac{1}{2} +l \right ]$.
According to \cite[Satz 3, p. 736]{Fi1966} we can take in (\ref{eq-3-1}) the $k$-th derivative with respect to $x$, which leads to
$$
	\int_{-\infty }^{\infty} {\rm{e}}^{nux} u^k \frac{1}{(1+{\rm{e}}^u)^n} g(u) du = 0 , \,
	 x \in \left [ \frac{1}{2}-l, \frac{1}{2} +l \right ] .
$$
For $x = \frac{1}{2}$ we get
\begin{equation}
\label{eq-3-2}
	\int_{-\infty }^{\infty} g(u) {\rm{e}}^{-\frac{|u|}{4}} u^k \frac{{\rm e}^{\frac{nu}{2}}}{(1+{\rm{e}}^u)^n}  
	 {\rm{e}}^{\frac{|u|}{4}}  du = 0 , \,
	k \in \mathbb{N}_0 .
\end{equation}
Let us remark that
\begin{equation}
\label{eq-3-3}
	\frac{{\rm e}^{\frac{nu}{2}}}{(1+{\rm{e}}^u)^n} {\rm{e}}^{\frac{|u|}{4}} \leq
	{\rm e}^{\frac{1-2n}{4}|u|} , \, u \in \mathbb{R} .
\end{equation}
According to \cite[Section 8.4.3, p. 428]{KoFo1975}, from (\ref{eq-3-3}) we deduce that the sequence
$$
	\left ( u^k   \frac{{\rm e}^{\frac{nu}{2}}}{(1+{\rm{e}}^u)^n} {\rm{e}}^{\frac{|u|}{4}}     \right )_{k \in \mathbb{N}_0}
$$
is complete in $ L^2 (\mathbb{R} )$.
Since $g \in C_b  (\mathbb{R} ) $, we have $g(u)  {\rm{e}}^{-\frac{|u|}{4}} \in  L^2 (\mathbb{R} ) $, and now 
(\ref{eq-3-2}) implies $g = 0$ a.~e. on $\mathbb{R}$. By using again the continuity of $g$ we get 
$g = 0$ on $\mathbb{R}$, and so $ f = 0$ on $[0,1]$.

\section{The eigenstructure of $\bar{\mathbb{B}}_n$}
\label{eigen}
By direct computation it is easy to find the first eigenvalues and eigenpolynomials of $\bar{\mathbb{B}}_n$:
\begin{align*}
	\eta_0^{(n)} & =  1 & q_0^{(n)}& = 1
	\\
	\eta_1^{(n)} & =  1 & q_1^{(n)}& = x-\frac{1}{2}
	\\
	\eta_2^{(n)} & =  \frac{n}{n+1} & q_2^{(n)}& = x(x-1)
	\\
	\eta_3^{(n)} & =  \frac{n^2}{(n+1)(n+2)} & q_3^{(n)}& = x(x-1) \left ( x-\frac{1}{2} \right )
	\\
	\eta_4^{(n)} & =  \frac{n^3}{(n+1)(n+2)(n+3)} & q_4^{(n)}& = x(x-1) \left ( x(x-1)+\frac{n+1}{5n+6} \right )
\end{align*}
As
\begin{eqnarray}
\label{BMonom}
	\bar{\mathbb{B}}_n e_0 &=& e_0 ,
	\\
\nonumber
	\bar{\mathbb{B}}_n e_k (x) & =& \frac{nx(nx+1) \dots (nx+k-1)}{n(n+1) \dots (n+k-1)} , \, k \geq 1,
\end{eqnarray} 
following directly from the definition of $\bar{\mathbb{B}}_n$, we conclude that the eigenvalues of 
$\bar{\mathbb{B}}_n \, : \, \Pi_n \longrightarrow \Pi_n$ are the numbers
$$
	\eta_k^{(n)}  =  \frac{(n-1)!}{(n+k-1)!} n^k , \, k \geq 0 .
$$
Let us denote by $p_k^{(n)}$ the eigenpolynomials of $B_n$ (see \cite{CoWa2000}). Here are some examples (see \cite[(9.1)]{CoWa2000}).
\begin{eqnarray*}
	 p_0^{(n)}& = &1
	\\
	p_1^{(n)}& = &x-\frac{1}{2}
	\\
	p_2^{(n)}& = &x(x-1)
	\\
	p_3^{(n)}& = &x(x-1) \left ( x-\frac{1}{2} \right )
	\\
	p_4^{(n)}& = &x(x-1) \left ( x(x-1)+\frac{n-1}{5n-6} \right )
\end{eqnarray*}
 Thus we have
$$
	q_k^{(n)} = p_k^{(n)} , \, 0 \leq k \leq 3
$$
and
\begin{equation}
\label{eq-4-1}
	\lim_{n \to \infty } q_k^{(n)} (x) = \lim_{n \to \infty } p_k^{(n)} (x) , \, k =4 , 
\end{equation}
uniformly in $[0,1]$.
We shall show that the eigenstructure of $\bar{\mathbb{B}}_n$ is similar to that of $B_n$; in particular, that (\ref{eq-4-1})
holds for all $k \geq 0$. Since the polynomials 
$$
	\lim_{n \to \infty } p_k^{(n)} (x):= p_k^* (x) , \, k \geq 0 ,
$$
are completely described in \cite{CoWa2000}, we get the same information about $\lim_{n \to \infty } q_k^{(n)} (x) $.
Let $k \geq 2$ and $n \geq 1$. We want to determine $q_k^{(n)} \in \Pi_k$ such that 
\begin{equation}
\label{eq-4-2}
	\bar{\mathbb{B}}_n q_k^{(n)} = \eta_k^{(n)} q_k^{(n)} .
\end{equation}
We put $ q_k^{(n)} (x) = \sum_{j=0}^k a(n,k,j) x^j$, with $a(n,k,k)=1$. Hence
$$
	\bar{\mathbb{B}}_n ( q_k^{(n)};x) = \sum_{j=0}^k a(n,k,j) \bar{\mathbb{B}}_n (e_j;x) .
$$
With (\ref{BMonom}) we derive
\begin{eqnarray}
\label{eq-4-3}
	\bar{\mathbb{B}}_n ( q_k^{(n)};x) 
	& = &
	\sum_{j=0}^k a(n,k,j) \frac{nx(nx+1) \dots (nx+j-1)}{n(n+1) \dots (n+j-1)}
	\\
	\nonumber
	&=&
	\frac{n^k}{n(n+1) \dots (n+k-1)} \sum_{j=0}^k a(n,k,j) x^j .
\end{eqnarray}
From the definition of the Stirling numbers of first kind $s(j,i)$, we obtain immediately
$$
	nx(nx+1) \dots (nx+j-1) = \sum_{i=0}^j s(j,i)(-1)^{j-i} n^i x^i ,
$$
so that (\ref{eq-4-3}) becomes, after some manipulation,
$$
	\sum_{i=0}^k \left \{ \sum_{j=i}^k \frac{ s(j,i)(-1)^{j-i} n^i}{n(n+1) \dots (n+j-1)}  a(n,k,j) \right \} x^i
	=
	\sum_{i=0}^k  \frac{a(n,k,i)  n^k }{n(n+1) \dots (n+k-1)}  x^i .
$$
This leads to
\begin{equation}
\label{eq-4-4}
	\sum_{j=i}^k \frac{ s(j,i)(-1)^{j-i} }{n(n+1) \dots (n+j-1)}  a(n,k,j) 
	=
	\frac{  n^{k-i} }{n(n+1) \dots (n+k-1)} a(n,k,i)  ,
\end{equation}
for all $i=0,1, \dots , k$.
Since $s(i,i) = 1$, we can solve (\ref{eq-4-4}) for $a(n,k,i)$ getting
\begin{equation}
\label{eq-4-5}
	a(n,k,i) = \frac{\sum_{j=i+1}^k (-1)^{j-i-1} s(j,i) (n+j)(n+j+1) \dots (n+k-1) a(n,k,j)}{ (n+i)(n+i+1) \dots (n+k-1) - n^{k-i}} ,
\end{equation}
for all $i \in \{k-1, k-2 , \dots , 0\}$.
Recalling that $n$ and $k$ are given, and $a(n,k,k) = 1$, (\ref{eq-4-5}) represents a recurrence relation for computing $a(n,k,i) $, $i = k-1, k-2, \dots , 0$. In particular, using $s(k,k-1) = -\frac{k(k-1)}{2}$, 
$s(k,k-2) = \frac{k(k-1)(k-2)(3k-1)}{24}$, we get
\begin{eqnarray}
\label{eq-4-6}
	a(n,k,k-1) 
	& = &
	- \frac{k}{2} ,
	\\
\label{eq4-7}
	a(n,k,k-2) 
	& = &
	\frac{k(k-1)(k-2)}{24} \cdot \frac{6n+3k-5}{(2k-3)n+(k-1)(k-2)}  .
\end{eqnarray}
Let us prove by induction that
\begin{equation}
\label{eq-4-8}
	a^* (k,j) := \lim_{n \to \infty} a(n,k,j) = \prod_{l=1}^{k-j} \frac{(k+1-l)(k-l)}{l(l-2k+1)} .
\end{equation}
For $j=k$ (\ref{eq-4-8}) is verified because $a(n,k,k)=1$.
Due to (\ref{eq-4-6}), (\ref{eq-4-8}) is verified also for $j=k-1$.
Suppose now that (\ref{eq-4-8}) is true for $j=i+1$, and let's prove it for $j=i$. From (\ref{eq-4-5}) we infer
\begin{eqnarray*}
	a(n,k,i)
	& = &
	\left \{ (i+(i+1) + \dots + (k-1))n^{k-i-1} + \mbox{ terms of lower degree} \right \}^{-1}
	\\
	& &
	\times s(i+1,i) \left (n^{k-i-1} + \mbox{ terms of lower degree} \right ) a(n,k,i+1) ,
\end{eqnarray*}  
so that, by the induction hypothesis,
\begin{eqnarray*}
	a^* (k,i)
	& = &
	\frac{s(i+1,i)}{i+(i+1)+ \dots +(k-1) } a^* (k,i+1)
	\\
	& = &
	- \frac{i(i+1)}{(k-i)(k+i-1)} \prod_{l=1}^{k-i-1} \frac{(k+1-l)(k-l)}{l(l-2k+1)}
	\\
	& = &
	 \prod_{l=1}^{k-i} \frac{(k+1-l)(k-l)}{l(l-2k+1)} ,
\end{eqnarray*}
and this completes the proof of (\ref{eq-4-8}).

It follows that
$$
	\lim_{n \to \infty} q_k^{(n)} (x) = \sum_{j=0}^{k} a^* (k,j) x^j ,
$$
and the coefficients $a^* (k,j)$ are equal to the coefficients $c^* (j,k)$ from \cite[Theorem 4.1]{CoWa2000}.
This leads to
\begin{equation}
\label{eq-4-9}
	\lim_{n \to \infty}  q_k^{(n)} (x) = \lim_{n \to \infty}  p_k^{(n)} (x) =: p_k^* (x) , \, k \geq 0 ,
\end{equation}
where(see \cite[Theorem 4.5]{CoWa2000}) $p_0^*(x) = 1$, $p_1^* (x) = x - \frac{1}{2}$, and 
\begin{equation}
\label{eq-4-10}
p_k^* (x) = \frac{k! (k-2)!}{(2k-2)!} x(x-1) P_{k-2}^{(1,1)} (2x-1) , \, k \geq 2 .
\end{equation}
($P_m^{(1,1)}$ are the Jacobi polynomials, orthogonal with respect to the weight $(1-t)(1+t)$ on the interval $[-1,1]$.)

\section{The operators $F_n$}
\label{operators}
The images of the monomials under $ \bar{\mathbb{B}}_n$ (see (\ref{BMonom})) show that 
$\bar{\mathbb{B}}_n \, : \, \Pi_n \longrightarrow \Pi_n$ is bijective.
 By composing the operators 
$$
	 B_n \, : \, C[0,1] \longrightarrow \Pi_n \mbox{ and } \bar{\mathbb{B}}_n^{-1} \, : \, \Pi_n \longrightarrow \Pi_n
$$
we obtain the operators
$$
	F_n :=  \bar{\mathbb{B}}_n^{-1} \circ B_n , \, F_n  \, : \, C[0,1] \longrightarrow \Pi_n , \, n \geq 1.
$$
Now $B_n$ can be represented as
$$
	B_n = \bar{\mathbb{B}}_n  \circ F_n, \, n \geq 1 .
$$
The eigenvalues of $B_n$ (see \cite{CoWa2000}) are
$$
	 \lambda_{k}^{(n)}= \frac{n!}{(n-k)!} \cdot \frac{1}{n^k} , \,0 \leq  k \leq n .
$$
It follows that the eigenvalues of  $F_n \, : \, \Pi_n \longrightarrow \Pi_n$ are
$$
	\nu^{(n)}_0
	 =\nu^{(n)}_1=1 , \, 
	 \nu^{(n)}_k=\frac{\lambda_{k}^{(n)}}{\eta^{(n)}_k}
	 =\frac{(n^2-1)(n^2-4) \dots (n^2-(k-1)^2)}{n^{2k-2}} , \,2 \leq  k \leq n ,
$$
or
$$
	 \nu^{(n)}_k=\frac{(n-1+k)!}{(n-k)!} \cdot \frac{1}{n^{2k-1}} , \,0 \leq  k \leq n.
$$
Here are some images of the monomials:
\begin{eqnarray*}
	F_n e_0 
	& = & 
	e_0 
	\\
	 F_n e_1  
	& = & 
	 e_1
	\\
	 F_n (e_2; x)
	& = & 
	\frac{x}{n^2} \left \{ (n^2-1)x+1 \right \}
	\\
	 F_n (e_3 ;x) 
	& = & 
	\frac{x}{n^4} \left \{ (n^2-1)(n^2-4)x^2+6(n^2-1)x+(2-n^2) \right \}
	\\
	 F_n (e_4; x) 
	& = & 
	\frac{x}{n^6} \left \{ (n^2-1)(n^2-4)(n^2-9)x^3 +18(n^2-1)(n^2-4)x^2 \right .
	\\
	& &
 	\qquad \left . -(n^2-1)(4n^2-n-42)x-(n^3+5n^2-n-6) \right \}
	\\
	 F_n (e_5 ;x) 
	& = & 
	\frac{x}{n^8} \left \{ (n^2-1)(n^2-4)(n^2-9)(n^2-16)x^4+40(n^2-1)(n^2-4)(n^2-9)x^3  \right .
	\\
	& &
 	\qquad  -5(n^2-1)(n^2-4)(2n^2-n-60)x^2-5(n^2-1)(n^3+13n^2-6n-72)x
	\\
	& &
 	\qquad \left . +2n^4-10n^3-25n^2+10n+24 \right \}
	\\
	 F_n (e_6; x) 
	& = &
	\frac{x}{n^{10}} \left \{  (n^2-1)(n^2-4)(n^2-9)(n^2-16)(n^2-25)x^5 \right .
	\\
	& &
 	\qquad +75(n^2-1)(n^2-4)(n^2-9)(n^2-16)x^4
	\\
	& &
 	\qquad -5(n^2-1)(n^2-4)(n^2-9)(4n^2-3n-260)x^3
	\\
	& &
 	\qquad 
	-15(n^2-1)(n^2-4)(n^3+25n^2-18n-360)x^2
	\\
	& &
 	\qquad 
	+(n^2-1)(22n^4-144n^3-919n^2+626n+3720)x
	\\
	& &
 	\qquad \left .+9n^5+16n^4-95n^3-135n^2+86n+120 \right \}
\end{eqnarray*}
In particular, from the representation of $F_n (e_3 ;x)$, $0 \leq x \leq 1$, we see that for $n \geq 2$ $F_n$ is {\bf not} a positive operator. Indeed we have
$$
	F_n \left (e_3 ; \frac{1}{(n+1)^2} \right ) =
	\frac{-n^5-6n^4-3n^3+14n^2+17n+6}{n^4(n+2)^5} < 0 , \, n \geq 2 .
$$
\begin{remark}
In Theorem \ref{thm1} we showed that $\bar{\mathbb{B}}_n \, : \, C[0,1] \longrightarrow C[0,1]$ is injective, hence 
$\bar{\mathbb{B}}_n^{-1}$ exists on the range $R(\bar{\mathbb{B}}_n)$.
If we assume that there is an operator $ Q \, : \, C[0,1] \longrightarrow C[0,1]$ such that 
$B_n f =( \bar{\mathbb{B}}_n \circ Q) f $ for all $f \in C[0,1]$, then we have
\begin{eqnarray*}
	F_n f 
	& = &
	( \bar{\mathbb{B}}_n^{-1} \circ B_n) f
	\\
	& = &
	[ \bar{\mathbb{B}}_n^{-1} \circ (  \bar{\mathbb{B}}_n \circ Q)] f
	\\
	& = &
	[ (\bar{\mathbb{B}}_n^{-1} \circ   \bar{\mathbb{B}}_n ) \circ Q] f
	\\
	& = &
	Qf ,
\end{eqnarray*}
so $F_n = Q$.
Since $F_n$ is not a positive operator, the equality shows that there is \underline{\bf no} positive operator allowing the decomposition in question.
\end{remark}

\section{The moments of $F_n$}
\label{moments}
Consider the moments of $F_n$, defined by 
$$
	M_{n,m}(x):= F_n((e_1-xe_0)^m ;x) , \, m \geq 0, \, x \in [0,1] .
$$
By using the above images of the monomials we get
\begin{eqnarray*}
	M_{n,0} (x)
	& = & 
	1
	\\
	 M_{n,1} (x)  
	& = & 
	 0
	\\
	 M_{n,2} (x)
	& = & 
	\frac{x(1-x)}{n^2}
	\\
	 M_{n,3} (x) 
	& = & 
	\frac{x(1-x)(1-2x)}{n^4}  \left \{ -n^2+2 \right \}
	\\
	 M_{n,4} (x) 
	& = & 
	\frac{x(1-x)}{n^6} \left \{3x(1-x) (11n^2-12)-n^3-5n^2+n+6 \right \}
	\\
	 M_{n,5} (x) 
	& = & 	
	\frac{x(1-x)(1-2x)}{n^8}\left \{-2x(1-x)(17n^4-160n^2+144) \right .
	\\
	& &
	\qquad \left . +2n^4-10n^3-25n^2+10n+24 \right \}
	\\
	 M_{n,6} (x) 
	& = &
	\frac{x(1-x)}{n^{10}} \left \{-5x^2(1-x)^2(8n^6-653n^4+3524n^2-2880) \right .
	\\
	& &
	\qquad 
	+5x(1-x)(2n^6-15n^5-155n^4+123n^3+872n^2-108n-720)
	\\
	& &
	\qquad \left . 
	+9n^5+16n^4-95n^3-135n^2+86n+120 \right \}
\end{eqnarray*}
In particular, 
\begin{eqnarray*}
	\lim_{n \to \infty} n^2 M_{n,2} (x) 
	& = &
	  x(1-x) , \,  \lim_{n \to \infty} n^2 M_{n,3} (x) =-x(1-x)(1-2x) , 
	\\
	\lim_{n \to \infty} n^2 M_{n,4} (x) 
	& = &
	0   .
\end{eqnarray*}
These facts, combined with Taylor's formula, lead to the following conjecture concerning a Voronovskaya-typ result.
\begin{conjecture}
\label{conj1}
For $f \in C^3[0,1]$ we have
$$
	\lim_{n\to \infty} n^2 (F_nf-f)(x) = \frac{x(1-x)}{2}f''(x)-\frac{x(1-x)(1-2x)}{6}f'''(x) ,
$$
uniformly on $[0,1]$.
\end{conjecture}
We will see that this conjecture is verified for all polynomials.

\section{A representation of $\bar{\mathbb{B}}_n^{-1} e_j$}
\label{repr1}
 Let us denote by $S(j,k)$ the Stirling numbers of second kind. From their definition we have
$$
	(-nx)^j = \sum_{k=0}^j S(j,k) (-nx)(-nx-1) \dots (-nx-k+1) ,
$$
i.~e., by using (\ref{BMonom})
$$
	(-1)^j n^j x^j =  \sum_{k=0}^j S(j,k)(-1)^k \bar{\mathbb{B}}_n (e_k ;x) \frac{(n-1+k)!}{(n-1)!} . 
$$ 
This entails
$$
	e_j = \frac{1}{n^j} \sum_{k=0}^j (-1)^{j-k} \frac{(n-1+k)!}{(n-1)!} S(j,k) \bar{\mathbb{B}}_n e_k  . 
$$ 
Finally we get
\begin{equation}
\label{eqfin}
	 \bar{\mathbb{B}}_n^{-1} e_j = \frac{1}{n^j} \sum_{k=0}^j (-1)^{j-k} \frac{(n-1+k)!}{(n-1)!} S(j,k) e_k   . 
\end{equation}
Here are some examples.
\begin{eqnarray*}
	 \bar{\mathbb{B}}_n^{-1} e_0
	& = &
	1 ,
	\\
	 \bar{\mathbb{B}}_n^{-1} e_1
	& = &
	e_1
	\\
	 \bar{\mathbb{B}}_n^{-1} e_2
	& = &
	\frac{n+1}{n} e_2 - \frac{1}{n} e_1
	\\
	 \bar{\mathbb{B}}_n^{-1} e_3
	& = &
	\frac{(n+1)(n+2)}{n^2}e_3-3\frac{n+1}{n^2} e_2 + \frac{1}{n^2} e_1
	\\
	 \bar{\mathbb{B}}_n^{-1} e_4
	& = &
	\frac{(n+1)(n+2)(n+3)}{n^3}e_4-6 \frac{(n+1)(n+2)}{n^3}e_3+7\frac{n+1}{n^3} e_2 - \frac{1}{n^3} e_1
\end{eqnarray*}
Since $S(j,j)=1$, it follows from (\ref{eqfin}) that
$$
	\lim_{n \to \infty}  \bar{\mathbb{B}}_n^{-1} e_j  =e_j , \, j \geq 0.
$$

\section{A first representation of $F_nf$}
\label{repr}
For $f \in C[0,1]$ we have
\begin{eqnarray*}
	F_nf
	 = 
	 \bar{\mathbb{B}}_n^{-1}(B_nf)
	& = &
	 \bar{\mathbb{B}}_n^{-1} \left ( \sum_{i=0}^n f\left(\frac{i}{n} \right ) b_{n,i} \right )
	\\
	& = &
	 \sum_{i=0}^n f\left(\frac{i}{n} \right ) \bar{\mathbb{B}}_n^{-1} b_{n,i} .
\end{eqnarray*}
Consider the polynomials $\varphi_{n,i}:= \bar{\mathbb{B}}_n^{-1} b_{n,i} $, $0 \leq i \leq n$. Then
$$
	F_nf =  \sum_{i=0}^n f\left(\frac{i}{n} \right ) \varphi_{n,i} , \, n \geq 1, \, f \in C[0,1].
$$
In fact,
\begin{eqnarray*}
	\varphi_{n,i}
	& = & 
	 \bar{\mathbb{B}}_n^{-1} \left ( \sum_{l=0}^{n-i}{n \choose i}{n-i \choose l} (-1)^l e_{i+l}  \right )
	\\
	& = &
	 \sum_{l=0}^{n-i}{n \choose i}{n-i \choose l} (-1)^l  \bar{\mathbb{B}}_n^{-1} \left ( e_{i+l}  \right ) .
\end{eqnarray*}
So we get
$$
	\varphi_{n,i}
	= 
	 \sum_{l=0}^{n-i} \sum_{k=0}^{i+l}{n \choose i}{n-i \choose l} \frac{(-1)^{i-k}}{n^{i+l}}
	n(n+1) \dots (n+k-1) S(i+l,k) e_k .
$$
Before giving some examples, we prove that 
\begin{equation}
\label{eq-8-1}
	\varphi_{n,i} (x) = \varphi_{n,n-i} (1-x) , ,\, 0 \leq i \leq n, \, x \in [0,1].
\end{equation}
Let $S \, : \, C[0,1] \longrightarrow C[0,1]$, $Sf(x) = f(1-x)$, $f \in C[0,1]$, $x \in [0,1]$.
\\
It is easy to see that
$$
	\bar{\mathbb{B}}_n \circ S = S \circ \bar{\mathbb{B}}_n .
$$	
It follows that $S = \bar{\mathbb{B}}_n^{-1} \circ S \circ \bar{\mathbb{B}}_n $, i.e.,
$$
	S \circ  \bar{\mathbb{B}}_n^{-1} =  \bar{\mathbb{B}}_n^{-1} \circ S .
$$
On the other hand, $b_{n,i}= S b_{n,n-i}$. Now
$$
	\varphi_{n,i} =  \bar{\mathbb{B}}_n^{-1} b_{n,i} 
	=  \bar{\mathbb{B}}_n^{-1} S b_{n,n-i} 
	= S \bar{\mathbb{B}}_n^{-1} b_{n,n-i} 
	= S  \varphi_{n,n-i} ,
$$
i.e., (\ref{eq-8-1}).

Here are some examples:
\begin{eqnarray*}
	\varphi_{1, 0}(x) 
	& = & 
	1-x
	\\
	\varphi_{2, 0}(x) 
	& = &
	(1-x)\left (1-\frac{3}{2}x \right )
	\\
 	\varphi_{2, 1}(x) 
	& = &
	3 x(1 -  x) 
	\\
	\varphi_{3, 0}(x) 
	& = &
	(1-x)(1-2x)\left ( 1-\frac{10}{9}x \right ) 
	\\
	\varphi_{3, 1}(x) 
	& = &
	 \frac{20}{3} x (1-x)  \left ( \frac{4}{5}-x \right )
	\\
	\varphi_{4, 0}(x)
	& = &
	(1-x) \left ( -\frac{105}{32} x^3 + \frac{225}{32}x^2 - \frac{305}{64} x +1 \right )
	\\
	\varphi_{4, 1}(x)
	& = &
	x(1-x) \frac{5}{16} (42x^2-66x+25)
	\\
	\varphi_{4, 2}(x)
	& = &
	x(1-x) \frac{15}{32} [42x(1-x)-5]
	\\
	\varphi_{5, 0}(x)
	& = &
	(1-x) \left ( \frac{3024}{625}x^4 - \frac{8736}{625}x^3+ \frac{9114}{625}x^2- \frac{4026}{625}+1 \right )
	\\
	\varphi_{5, 1}(x)
	& = &
	x(1-x)  \frac{6}{125}(-504x^3+1176x^2-889x+216)
	\\
	\varphi_{5, 2}(x)
	& = &
	x(1-x)  \frac{24}{125} (252x^3-448x^2+217x-18)
\end{eqnarray*}
\begin{eqnarray*}
	\varphi_{6, 0}(x)
	& = &
	(1-x)  \frac{1}{108}\left (-770x^5+2800x^4-3955x^3+2695x^2- \frac{63217}{72}x + 108 \right )
	\\
	\varphi_{6, 1}(x)
	& = &
	x(1-x)  \frac{7}{54} \left ( 330x^4-1020x^3+1155x^2-565x+ \frac{2401}{24} \right )
	\\
	\varphi_{6, 2}(x) 
	& = &
	x(1-x)  \frac{35}{108} \left (- 330x^4+840x^3-723x^2+227x- \frac{343}{24} \right )
	\\
	\varphi_{6, 3}(x)
	& = &
	x(1-x)  \frac{35}{27} \left [x(1-x) (110x(1-x)-23)+ \frac{49}{72} \right ]
	\\
\end{eqnarray*}
The  \textquotedblright Lebesgue function\textquotedblright corresponding to $F_n$ is
$$
	\Psi_n (x) = \sum_{i=0}^n | \varphi_{n,i} (x)|, \, x \in [0,1].
$$

\begin{figure}[h]
\includegraphics[scale=0.6]{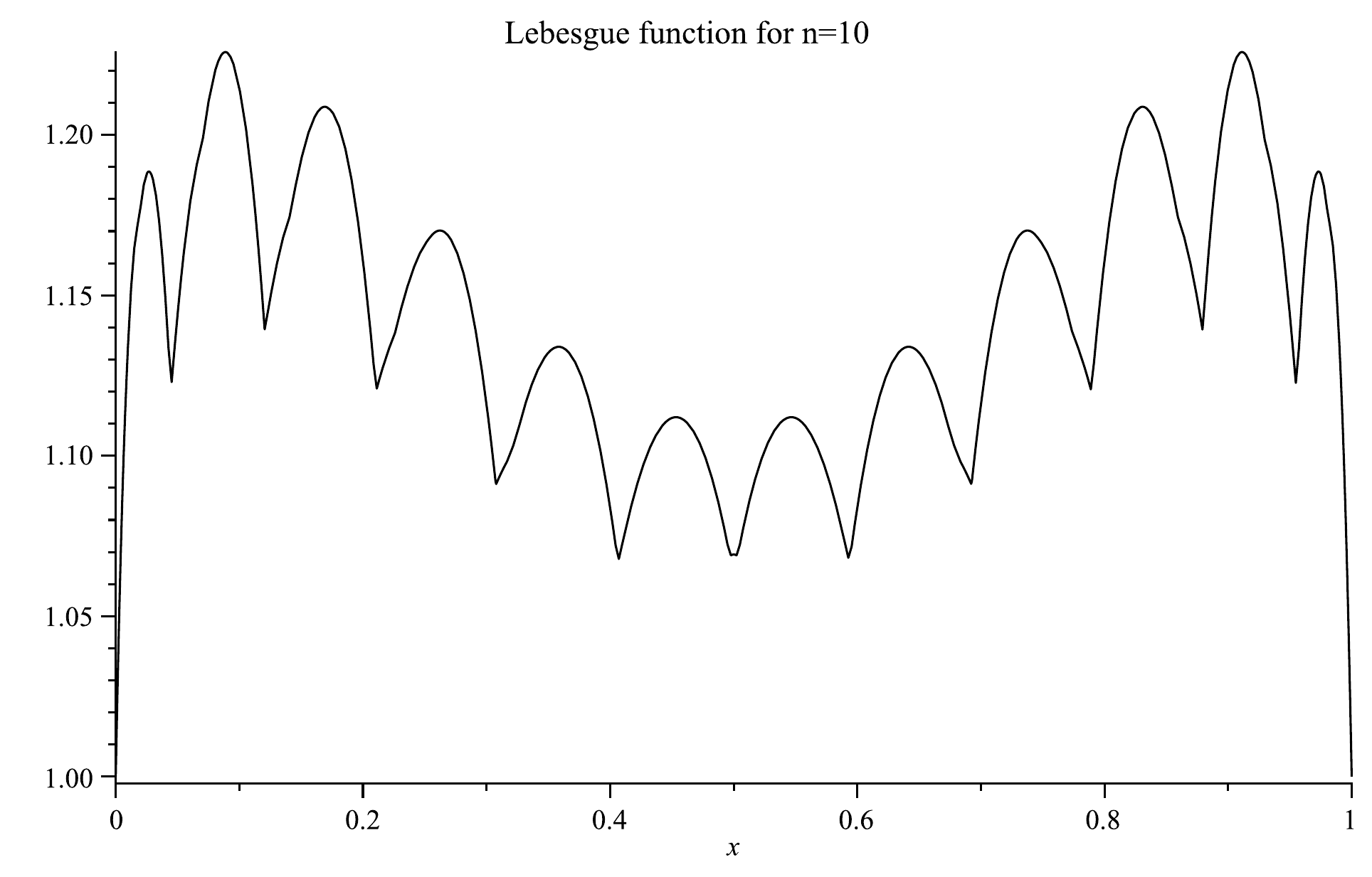}
\end{figure}
\begin{figure}[h]
\includegraphics[scale=0.6]{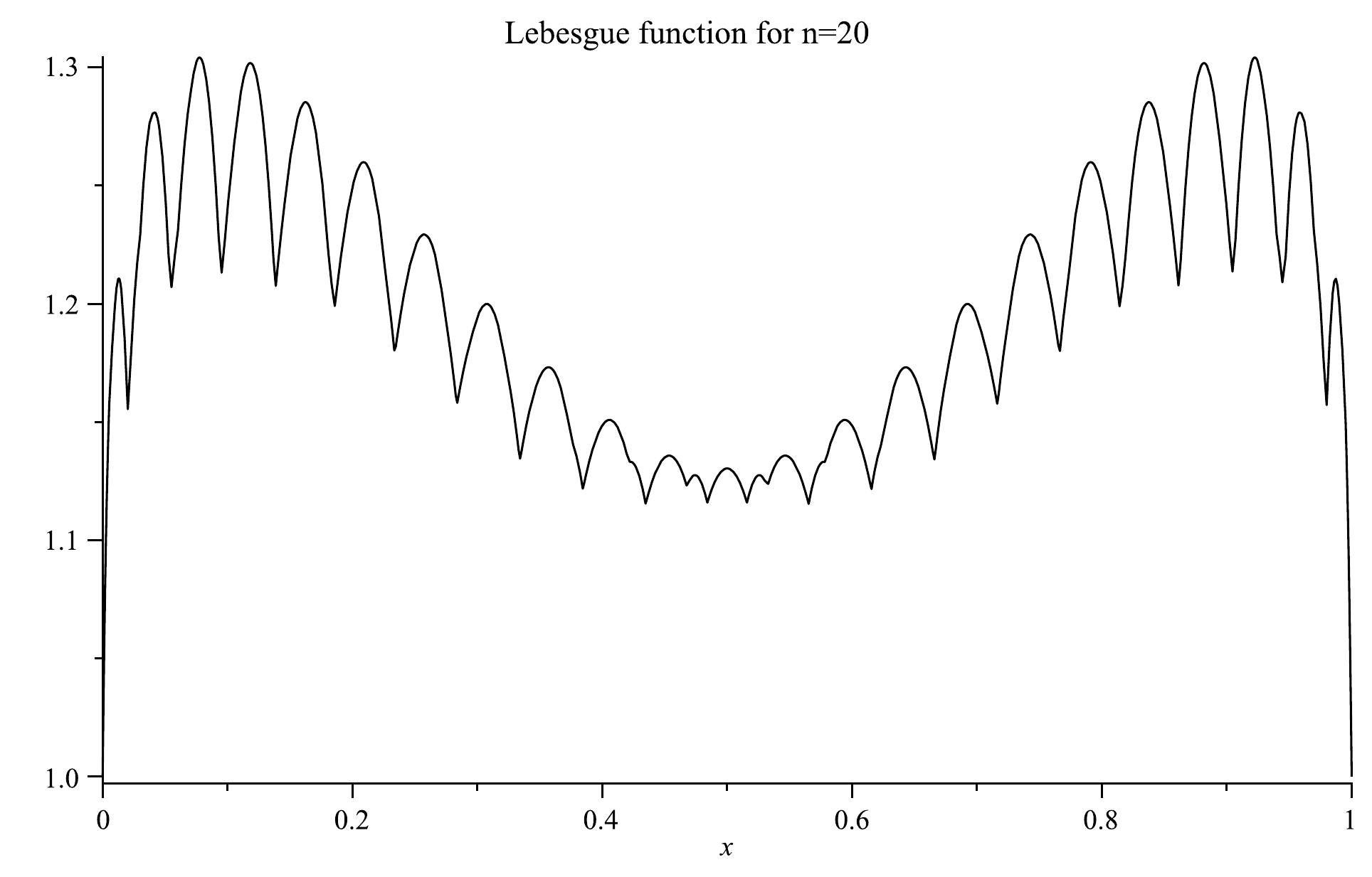}
\end{figure}
\begin{figure}[h]
\includegraphics[scale=0.6]{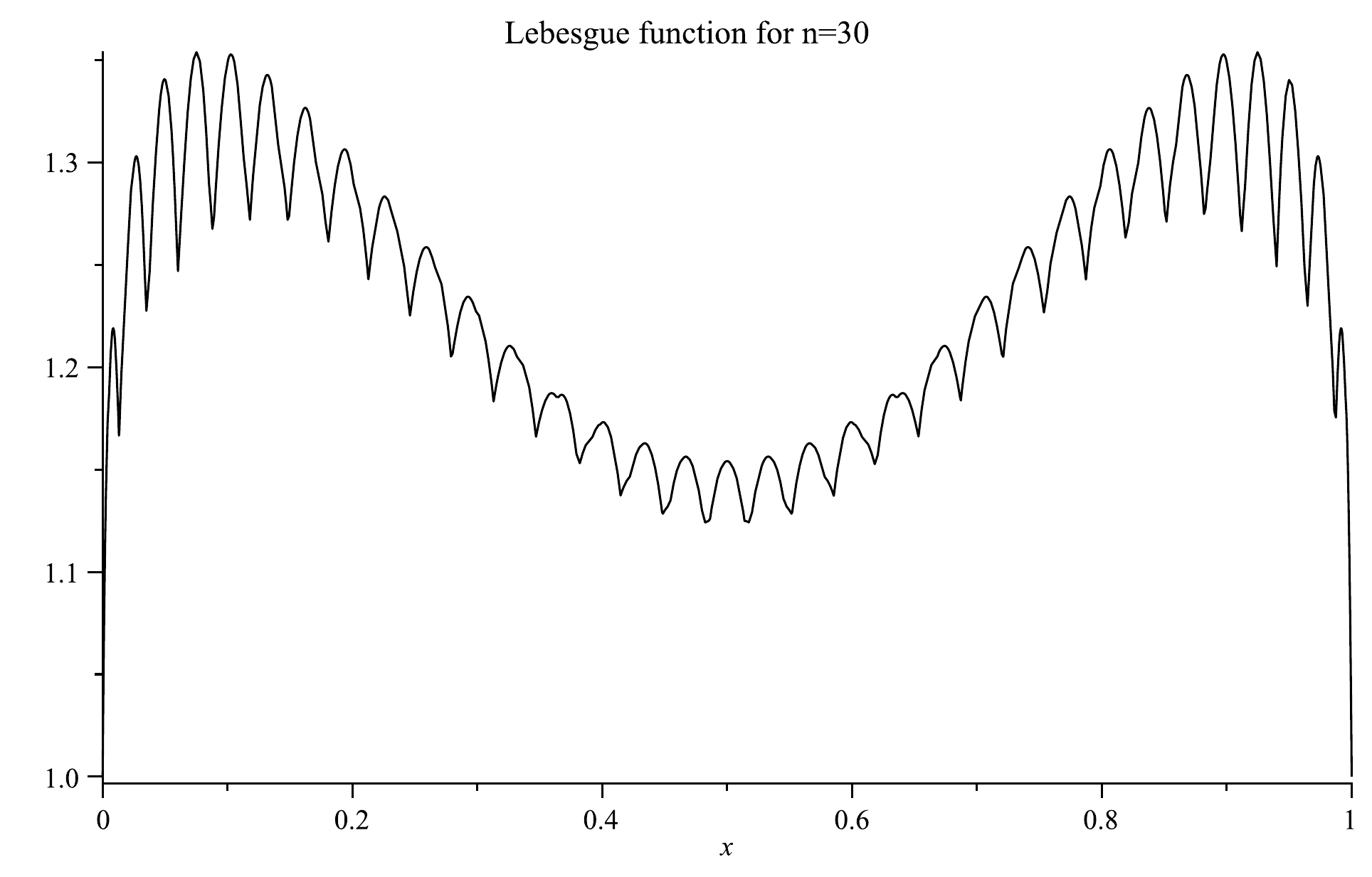}
\end{figure}

Experimental maximum values of the Lebesgue function for different $n$.

\begin{tabular}{c||c|c|c|c|c|c|c|}
	n & $10$ & $20$ & $30$ & $40$ & $50$ & $60$ & $70$
	\\ \hline
	Max. & $1.266$ & $1.304$ & $1.354$ & $1.387$ & $1.409$ & $1.433$ & $1.459$
\end{tabular}

\section{A second representation of $F_nf$}
\label{repr2}
Let us begin with a second representation of $B_nf$, namely
$$
	B_nf = \sum_{j=0}^n {n \choose j} \Delta_{\frac{1}{n}}^j f(0)  e_j .
$$
Now, with the results of Section \ref{repr1},
\begin{eqnarray}
\nonumber
	F_n f 
	& = &
	 \bar{\mathbb{B}}_n^{-1} (B_n f) 
	\\
\nonumber
	& = &
	\sum_{j=0}^n {n \choose j} \Delta_{\frac{1}{n}}^j f(0)  \bar{\mathbb{B}}_n^{-1} e_j
	\\
\nonumber
	& = &
	\sum_{j=0}^n {n \choose j}\frac{ j!}{ n^{j}} [0,\frac{1}{n}, \dots ,\frac{j}{n};f ]
	\sum_{k=0}^j  (-1)^{j-k} \frac{(n-k+1)!}{(n-1)!} \frac{1}{n^j} S(j,k)  e_k   
	\\
\label{eq16a}
	& = &
	\sum_{j=0}^n  [0,\frac{1}{n}, \dots ,\frac{j}{n};f ]  \frac{1}{(n-j)!  n^{2j-1}}
	\sum_{k=0}^j  (-1)^{j-k} (n+k-1)!  S(j,k)  e_k .
\end{eqnarray}
Consider the polynomials
$$
	\rho_{n,j}:=  \frac{1}{(n-j)!  n^{2j-1}}
	\sum_{k=0}^j  (-1)^{j-k} (n+k-1)!  S(j,k)  e_k   , \, 0 \leq j \leq n.
$$
Then
$$
	F_nf = \sum_{j=0}^n  [0,\frac{1}{n}, \dots ,\frac{j}{n};f ] \rho_{n,j} , \, f \in C[0,1] .
$$
Here are some examples concerning $\rho_{n,j}$.
\begin{eqnarray*}
	\rho_{n,0} (x)
	& = & 
	1
	\\
	 \rho_{n,1}  (x)
	& = & 
	 x
	\\
	 \rho_{n,2} (x)
	& = & 
	\frac{n-1}{n^2} x \left [ (n+1)x-1 \right ]
	\\
	 \rho_{n,3} (x) 
	& = & 
	\frac{(n-1)(n-2)}{n^4} x \left [ (n+1)(n+2) x^2 -3(n+1) x +1 \right ]
	\\
	 \rho_{n,4} (x) 
	& = & 
	\frac{(n-1)(n-2)(n-3)}{n^6} x \left [ (n+1)(n+2)(n+3)x^3 \right .
	\\
	& &
	\qquad \qquad \left . -6(n+1)(n+2)x^2 +7(n+1)x-1 \right ] .
\end{eqnarray*}

\section{An asymptotic formula for $F_n p$, $p \in \Pi$}
\label{asmp}
It is known that $S(m,j) = 0$ for $j >m$, and 
$  [0,\frac{1}{n}, \dots ,\frac{j}{n};e_m ]=n^{j-m} S(m,j) $.
Consequently, from (\ref{eq16a}) and
$S(m,m-1) = \frac{1}{2}m(m-1)$, $S(m,m-2)= \frac{1}{24}m(m-1)(m-2)(3m-5)$ we get
\begin{eqnarray*}
	F_n( e_m ;x) 
	& = &
	\sum_{j=0}^m
	\sum_{k=0}^j  (-1)^{j-k}  S(m,j) S(j,k) \frac{ (n+k-1)!}{(n-j)! n^{m+j-1}}  x^k   
	\\
	& = &
	\frac{1}{n^m} \sum_{k=0}^m x^k
	\sum_{j=k}^m  (-1)^{j-k}  S(m,j) S(j,k)  \left [ n^k - \frac{(j-k)(j+k-1)}{2}n^{k-1} \right .
	\\
	& &
	\qquad + \frac{1}{4}  \left  (\frac{(j-k)^2(j+k-1)^2}{2}-   \frac{(k-1)k(2k-1)+(j-1)j(2j-1)}{3} \right ) n^{k-2}
	\\
	& &
	\qquad  +\left . \mbox{ terms of degree } < k-2 \right ]
	\\
	& = &
	\frac{1}{n^m} \left \{ x^m \left [ n^m - \frac{m(m-1)(2m-1)}{6}n^{m-2}+ \mbox{ terms of degree } < m-2 \right ] \right .
	\\
	& &
	\qquad + x^{m-1} \left [\frac{m(m-1)^2}{2} n^{m-2}+ \mbox{ terms of degree } < m-2 \right ] 
	\\
	& &
	\qquad + x^{m-2} \left [-\frac{m(m-1)(m-2)}{6} n^{m-2}+ \mbox{ terms of degree } < m-2 \right ] 
	\\
	& &
	\qquad \left . +  \mbox{ terms of degree } < m-2 \right \} . 
\end{eqnarray*}
Thus 
\begin{eqnarray*}
	n^m (F_n( e_m ;x) -x^m) 
	& = &
	n^{m-2} \frac{m(m-1)}{6} x^{m-2} (1-x) ((2m-1)x-m+2) 
	\\
	& & + \mbox{ terms of degree } < m-2 .
\end{eqnarray*}
It follows that for each $m \geq 0$
$$
	\lim_{n\to \infty} n^2 (F_n(e_m;x)-x^m) = \frac{x(1-x)}{2}e''_m(x)-\frac{x(1-x)(1-2x)}{6}e'''_m(x) ,
$$
uniformly on $[0,1]$.
This implies that for any polynomial $p$ we have
$$
	\lim_{n\to \infty} n^2 (F_n(p;x)-p (x)) = \frac{x(1-x)}{2}p''(x)-\frac{x(1-x)(1-2x)}{6}p'''(x) .
$$
What is remarkable here is the factor $n^2$ (where an $n$ might have been expected).
Thus Conjecture \ref{conj1} is verified for all $p \in \Pi$.
In particular
$$
	\lim_{n\to \infty} F_n p =p , \, p \in \Pi .
$$
Moreover, by using \cite[(6.4)]{GoKaRa2007}, we get also             
$$
	\lim_{n\to \infty}  n^2 (F_n p- p) = 2 \lim_{n\to \infty}   n^2 (B_n p- U_{2n}p) , \, p \in \Pi .
$$

\section{Connection with a conjecture of Sh.~Cooper and Sh.~Waldron}
\label{conn}
Let 
$$
	p_0^{(n)} (x):=1, \, p_1^{(n)}(x):=x-\frac{1}{2}, \, p_2^{(n)}(x), \dots , p_n^{(n)}(x)
$$
be the monic eigenpolynomials of $B_n$, corresponding to the eigenvalues
$$
	\lambda_{0}^{(n)}=\lambda_{1}^{(n)}=1 , \, \lambda_{k}^{(n)}= \frac{(n-1) \dots (n-k+1)}{n^{k-1}} , \,2 \leq  k \leq n .
$$
Then 
$$
	B_nf = \sum_{k=0}^{n}  \lambda_{k}^{(n)}  p_k^{(n)} \mu_k^{(n)} (f) , \, f \in C[0,1] ,
$$
where $\left ( \mu_k^{(n)}  \right )_{0 \leq k \leq n}$ are the dual functionals (see \cite[Theorem 2.3]{CoWa2000}.
It is known that for each $k \geq 0$, 
\begin{eqnarray*}
	\lim_{n \to \infty} p_k^{(n)} 
	& = &
	p_k^* \in \Pi_k , \mbox{ uniformly on } [0,1] ,
	\\
	\lim_{n \to \infty} \mu_k^{(n)} (p) 
	& = &
	\mu_k^* (p) , \, p \in \Pi  ,
	\\
	 \lim_{n \to \infty} \lambda_k^{(n)}  
	& = &
	1
\end{eqnarray*}
(see \cite[Theorem 4.1, Theorem 4.20]{CoWa2000}). Moreover, according to  \cite[(4.18)]{CoWa2000},
$$
	p = \sum_{k=0}^{s} p_k^* \mu_k^* (p) , \, p \in \Pi_s .
$$
In particular,
$$
	p_s^* = \sum_{k=0}^{s} p_k^* \mu_k^* (p_s^*) , 
$$
and from the linear independence of $p_0^* , \dots , p_s^*$ we derive
\begin{equation}
\label{eq-11-1}
	\mu_s^* ( p_s^*) = 1 , \, s \geq 0 .
\end{equation}
Now 
$$
	B_n  p_s^*  =  \sum_{k=0}^{s} \lambda_{k}^{(n)} p_k^{(n)} \mu_k^{(n)} (p_s^*) 
$$
and 
$$
	F_n  p_s^*  =  \sum_{k=0}^{s} \lambda_{k}^{(n)} \mu_k^{(n)} (p_s^*)  \bar{\mathbb{B}}_n^{-1} ( p_k^{(n)}) .
$$
Since 
$$
	\lim_{n \to \infty} F_n  p_s^* =  p_s^* ,
$$
we conclude that
\begin{equation}
\label{eq-11-2}
	\lim_{n \to \infty}  \sum_{k=0}^{s} \lambda_{k}^{(n)} \mu_k^{(n)} (p_s^*)  \bar{\mathbb{B}}_n^{-1} ( p_k^{(n)}) 
	=
	 \sum_{k=0}^{s} \mu_k^* (p_s^*)  p_k^* , \, s \geq 0 .
\end{equation}
We know that 
$$
	 \bar{\mathbb{B}}_n^{-1} (p_0^{(n)}) = p_0^* \mbox{ and }  \bar{\mathbb{B}}_n^{-1} (p_1^{(n)}) = p_1^*  .
$$
Writing (\ref{eq-11-1}) and  (\ref{eq-11-2}) for $s=2$, we get
$$
	\lim_{n \to \infty}  \bar{\mathbb{B}}_n^{-1} (p_2^{(n)}) = p_2^* .
$$ 
Using (\ref{eq-11-1}) and  (\ref{eq-11-2}) with $s=3, 4, \dots $, we obtain
\begin{equation}
\label{eq-11-3}
	\lim_{n \to \infty}  \bar{\mathbb{B}}_n^{-1} (p_s^{(n)}) = p_s^*  , \, s \geq 0 .
\end{equation}
Consequently,
$$
	\lim_{n \to \infty} F_n (p_s^{(n)}) =
	\lim_{n \to \infty}  \bar{\mathbb{B}}_n^{-1} (B_n (p_s^{(n)})) =	
	\lim_{n \to \infty}  \lambda_s^{(n)} \bar{\mathbb{B}}_n^{-1}  (p_s^{(n)}) =	
	 p_s^*  ,
$$ 
i.~e.,
\begin{equation}
\label{eq-11-4}
	\lim_{n \to \infty} F_n (p_s^{(n)}) = p_s^*  ,  \, s \geq 0 .
\end{equation}
Concerning (\ref{eq-11-3}) and  (\ref{eq-11-4}) see also Section \ref{asymp}. 
In \cite[Remark on p. 149]{CoWa2000} the authors conjecture that
$$
	\lim_{n \to \infty}  \mu_k^{(n)} (f) = \mu_k^* (f) , \, f \in C[0,1] .
$$
\underline{Suppose that this is true.} If for a function  $ f \in C[0,1]$, 
\begin{eqnarray*}
	\lim_{n \to \infty} \sum_{k=0}^{n} \lambda_{k}^{(n)} \mu_k^{(n)} (f)  \bar{\mathbb{B}}_n^{-1} ( p_k^{(n)})
	& = &
	 \sum_{k=0}^{\infty}  \lim_{n \to \infty} \lambda_{k}^{(n)} \mu_k^{(n)} (f)  \bar{\mathbb{B}}_n^{-1} ( p_k^{(n)}),
\end{eqnarray*}
then
\begin{eqnarray*}
	\lim_{n \to \infty} F_n f
	& = &
	 \sum_{k=0}^{\infty}  \mu_k^* (f) p_k^*
\end{eqnarray*}
In the setting of \cite[Lemma 4.10]{CoWa2000}, the last series represents the function $f$.

\section{The asymptotic behavior of  $\bar{\mathbb{B}}_n^{-1} p_k^{(n)}- p_k^{(n)}$ and \\$F_n p_k^{(n)}- p_k^{(n)}$}
\label{asymp}
According to \cite[Theorem 2.3]{CoWa2000}, the eigenpolynomials of $B_n$ are
$$
	p_k^{(n)} = \sum_{j=0}^k c(j,k,n)  e_j , \, 0 \leq k \leq n .
$$
Moreover, according to \cite[Theorem 4.1]{CoWa2000},
\begin{eqnarray*}
	\lim_{n \to \infty} c(j,k,n)
	= 
	c^* (j,k) &=& \prod_{i=1}^{k-j} \frac{(k+1-i)(k-i)}{i(i-2k+1)} \mbox{ if }
	(j,k) \not= (0,1),
	\\
	c^*(0,1) & = & -\frac{1}{2}.
\end{eqnarray*}
Let $k \geq 2$. With the results of Section \ref{repr1} we get
\begin{eqnarray*}
	& &
	\bar{\mathbb{B}}_n^{-1} p_k^{(n)}
	\\
	 & = &
	 \sum_{j=0}^k c(j,k,n) \bar{\mathbb{B}}_n^{-1} e_j
	\\
	& = &
	 \sum_{j=0}^k c(j,k,n) \frac{1}{n^j} 
	\sum_{i=0}^j (-1)^{j-i}\frac{(n-1+i)!}{(n-1)!} S(j,i) e_i
	\\
	& = &
	\frac{1}{n^k} \sum_{j=0}^k c(j,k,n) n^{k-j} 
	\left \{ [\frac{(n-1+j)!}{(n-1)!} e_j \right .
	\\
	& & 
	\qquad \qquad \qquad - \frac{(n-2+j)!}{(n-1)!} \frac{j(j-1)}{2} e_{j-1} ]
	\\
	& & 
	\qquad
	\left . + \sum_{i=0}^{j-2} (-1)^{j-i} \frac{(n-1+i)!}{(n-1)!} S(j,i) e_i \right \}
	\\
	& = &
	\frac{1}{n^k} \sum_{j=0}^k c(j,k,n) n^{k-j}  
	\left \{ n^j  e_j +n^{j-1}  \frac{j(j-1)}{2} e_j \right .
	\\
	& & 
	\qquad \qquad \qquad \left . - n^{j-1} \frac{j(j-1)}{2} e_{j-1}+ q_{j-2} (n) \right \} ;
\end{eqnarray*}
here, as a polynomial in $n$, $  q_{j-2}$ has degree $j-2$, $j \geq 2$.
\begin{eqnarray*}
	\bar{\mathbb{B}}_n^{-1} p_k^{(n)}
	& = & 
	\underbrace{\sum_{j=0}^k c(j,k,n)   e_j}_{= p_k^{(n)}} + \frac{1}{n}  \sum_{j=2}^k c(j,k,n)  \frac{j(j-1)}{2}  (e_j-e_{j-1}) 
	\\
	& & 
	 +\sum_{j=2}^k  \frac{1}{n^j}    q_{j-2} (n) .
\end{eqnarray*}
\begin{eqnarray*}
	\bar{\mathbb{B}}_n^{-1} p_k^{(n)} -  p_k^{(n)}
	& = &
	 \frac{1}{n}  \left (  \sum_{j=2}^k c(j,k,n)  \frac{j(j-1)}{2}  (e_j-e_{j-1}) 
	+ \sum_{j=2}^k  \frac{1}{n^{j-1}}    q_{j-2} (n) \right ) .
\end{eqnarray*}
We get finally,
$$
	\lim_{n \to \infty} n \left ( \bar{\mathbb{B}}_n^{-1} p_k^{(n)} -  p_k^{(n)}  \right )
	=
	  \sum_{j=2}^k c^*(j,k)  \frac{j(j-1)}{2}  (e_j-e_{j-1}) , \, k \geq 2 . 
$$
Since $p_0^{(n)}=e_0$ and $p_1^{(n)}=e_1-\frac{1}{2} e_0$, we have also
$$
	\lim_{n \to \infty} n \left ( \bar{\mathbb{B}}_n^{-1} p_k^{(n)} -  p_k^{(n)}  \right )
	=
	 0 , \, k =0, \, 1 . 
$$
This leads to 
$$
	\lim_{n \to \infty} n \left ( \bar{\mathbb{B}}_n^{-1} p_k^{(n)}(x) -  p_k^{(n)} (x)  \right )
	=
	  - \frac{x(1-x)}{2} \sum_{j=0}^k c^*(j,k)  j(j-1)  x^{j-2} , \, k \geq 0 , 
$$
uniformly on $[0,1]$.
Moreover,
\begin{eqnarray*}
	F_n  p_k^{(n)} -  p_k^{(n)}
	& = &
	\bar{\mathbb{B}}_n^{-1}( B_n p_k^{(n)}) -  p_k^{(n)} 
	\\
	& = & \lambda_k^{(n)}  \bar{\mathbb{B}}_n^{-1} p_k^{(n)} -  p_k^{(n)} ,
\end{eqnarray*}
\begin{eqnarray*}
	n (F_n  p_k^{(n)} -  p_k^{(n)})
	& = &
	n [  \lambda_k^{(n)}  ( \bar{\mathbb{B}}_n^{-1}  p_k^{(n)} -  p_k^{(n)} ) + 
	( \lambda_k^{(n)}-1)   p_k^{(n)} ] .
\end{eqnarray*}
Since
\begin{eqnarray*}
	\lim_{n \to \infty}  \lambda_k^{(n)} 
	& = &
	1  , \,  \lim_{n \to \infty} n ( \lambda_k^{(n)}-1) = - \frac{k(k-1)}{2} ,
	\\
	\lim_{n \to \infty}    p_k^{(n)} (x)
	& = &
	p_k^* (x) =  \sum_{j=0}^k c^*(j,k)   x^{j}
\end{eqnarray*}
we get
\begin{eqnarray*}
	& &
	\lim_{n \to \infty}  n ( F_n  p_k^{(n)} (x) -  p_k^{(n)} (x) )
	\\
	& = &
	\frac{1}{2}  \sum_{j=0}^k (j-1)j (x^j - x^{j-1}) c^*(j,k) 
	- \frac{1}{2}  \sum_{j=0}^k (k-1)k x^j  c^*(j,k) 
	\\
	& = &
	\frac{1}{2}  \sum_{j=0}^k [(j-1)j (x^j - x^{j-1})-k(k-1) x^j ] c^*(j,k) .
\end{eqnarray*}

\section{A different approach to $F_n$}
\label{approach}
The genuine Bernstein-Durrmeyer operator $U_n$ can be described as $U_n = B_n \circ \bar{\mathbb{B}}_n $.
It follows that
$U_n \circ \bar{\mathbb{B}}_n^{-1} =B_n$ and $U_n  \circ \bar{\mathbb{B}}_n^{-1} \circ B_n = B_n \circ B_n$.
This leads to $U_n \circ F_n = B_n^2$, i.e.,
$$
	F_n = U_n^{-1} \circ B_n^2 .
$$
We will show that the inverse of $U_n$ for $\Pi_n$ can be written as
\begin{equation}
\label{inveq1}
	U_n^{-1} p = \sum_{l=0}^{n-1} (-1)^l \frac{(n-1-l)!}{l!(n-1)!} \widetilde{D}^{2l} p , \, p \in \Pi_n ,  
\end{equation}
where $ \widetilde{D}^{0}= I$, $ \widetilde{D}^{2l} = D^{l-1} \left [ x^l (1-x)^l D^{l+1}\right ]$, $l \geq 1$.
\\
According to \cite[Theorem 4]{GoSh1991} the eigenpolynomials of $U_n$ are 
\begin{equation}
\label{inveq2}
	p_0 (x) = 1 , \, p_1 (x) = x , \, p_k (x) = D^{k-2} \left [ x^{k-1} (1-x)^{k-1} \right ] , \, k \geq 2 ,
\end{equation}
with corresponding eigenvalues
$$
	\omega_k^{(n)} = \frac{(n-1)! n!}{(n-k)! (n+k-1)!} , \, 0 \leq k \leq n .
$$
It is shown in \cite[Lemma 2]{HeWa2011} that the differential operators $\widetilde{D}^{2l}$ posses the same eigenpolynomials with corresponding eigenvalues
\begin{align*}
\gamma_{k}^{(l)}:=
\begin{cases}
(-1)^l\frac{(k-1+l)!}{(k-1-l)!} , & 0 \leq l \leq k-1 ,\\
0 , &  l > k-1 .
\end{cases}
\end{align*}
So, to prove (\ref{inveq1}), we have to show that
$$
	U_n^{-1} p_k = \frac{1}{\omega_k^{(n)}} p_k , 0 \leq k \leq n .
$$
Indeed, we have
\begin{eqnarray*}
	U_n^{-1} p_k
	& = &
	\sum_{l=0}^{n-1} (-1)^l   \frac{(n-1-l)!}{l!(n-1)!} \widetilde{D}^{2l} p_k 
	\\
	& = &
	p_k \sum_{l=0}^{n-1} (-1)^l   \frac{(n-1-l)!}{l!(n-1)!} \gamma_{k}^{(l)}
	\\
	& = &
	\frac{1}{\omega_k^{(n)}} p_k ,
\end{eqnarray*}
where the last equation follows from the proof of \cite[Theorem 1]{HeWa2011}.

From $U_n(F_nf) = B_n(B_nf) $ we get
\begin{eqnarray*}
	&&
	 F_n f(0) b_{n,0} + F_n f(1) b_{n,n} + (n-1) \sum_{k=1}^{n-1} b_{n,k} \int_0^1 b_{n-2,k-1} (t) F_n f (t) dt
	\\
	& = &
	B_n f(0) ) b_{n,0} + B_n f(1) b_{n,n} + \sum_{k=1}^{n-1} b_{n,k} (B_n f) \left ( \frac{k}{n} \right ) .
\end{eqnarray*}
Consequently,
$$
	(n-1)  \int_0^1 b_{n-2,k-1} (t) F_n f (t) dt =  (B_n f) \left ( \frac{k}{n} \right ) , \, k = 1, \dots , n-1 ,
$$
which entails also
$$
	  \int_0^1 (1-t^{n-2}-(1-t)^{n-2}) F_n f (t) dt = \frac{1}{n-1} \sum_{k=1}^{n-1}   (B_n f) \left ( \frac{k}{n} \right )  .
$$
On the other hand,
$$
	F_n f = \sum_{i=0}^n  f \left ( \frac{i}{n} \right ) \varphi_{n,i} ,
$$
and so
$$
	(n-1) \sum_{i=0}^{n}  \left ( \int_0^1 b_{n-2,k-1} (t) \varphi_{n,i} (t) dt \right ) f \left ( \frac{i}{n} \right )
	=
	\sum_{i=0}^{n}  b_{n,i}  \left ( \frac{k}{n} \right )  f \left ( \frac{i}{n} \right ) ,
$$
for all $f \in C[0,1]$, $k=1 , \dots , n-1$. This gives
$$
	 \int_0^1 b_{n-2,k-1} (t) \varphi_{n,i} (t) dt = \frac{1}{n-1}  b_{n,i} \left ( \frac{k}{n} \right ) ,
$$
for all $k=1 , \dots , n-1$ and $i=0 , 1,  \dots , n$.

Let $J_i (x) := J_i^{(1,1)} (x)$, $i=0,1, \dots $, be the Jacobi polynomials on $[0,1]$, characterized by
$$
	\int_0^1 J_i (x) J_j (x) x(1-x) dx = \delta_{ij} , \, i,j \geq 0 .
$$
For $k \geq 2$ we rewrite the eigenpolynomials of $U_n$ in (\ref{inveq2}) into $p_k (x) =x(1-x) J_{k-2} (x)$.
For each $f \in C[0,1]$ the polynomial $B_n f$ can be represented as
$$
	B_n f = \sum_{k=0}^n \omega_k^{(n)} \nu_k^{(n)} (f) p_k ,
$$
and this representation introduces the linear functionals $ \nu_k^{(n)} \, : \, C[0,1] \longrightarrow \mathbb{R}$,
$ k=0, \dots , n $.

We also have
$$
	B_n^2 f = B_n (B_nf) =  \sum_{k=0}^n \omega_k^{(n)} \nu_k^{(n)} (B_n f) p_k , \, f \in C[0,1] .
$$
Since $U_n^{-1}p_k = (\omega_k^{(n)}  )^{-1} p_k$, it follows that
$$
	F_n f = U_n^{-1} B_n^2 f =  \sum_{k=0}^n  \nu_k^{(n)} (B_n f) p_k , \, f \in C[0,1] .
$$
On the other hand, let 
$$
	L f(x) := (1-x) f(0) + x f(1) .
$$
It is easy to see that $ \nu_0^{(n)} (f) = f(0)$, $ \nu_1^{(n)} (f) = f(1)-f(0)$, so that
$$
	B_n f(t) = L f(t) +  \sum_{k=2}^n \omega_k^{(n)} \nu_k^{(n)} ( f) t (1-t) J_{k-2} (t) , \, f \in C[0,1] , \, t \in [0,1].
$$
This entails 
$$
	\int_0^1 (B_n f(t) -L f(t) ) J_{k-2} (t) dt = \omega_k^{(n)} \nu_k^{(n)} ( f)  , \, f \in C[0,1] , \,  k=2, \dots , n .
$$
So we have an explicit description of the functionals $  \nu_k^{(n)}  $:
\begin{eqnarray*}
	 \nu_0^{(n)} (f) 
	& = &
	 f(0), \,   \nu_1^{(n)} (f) = f(1)-f(0) , 
	\\
	 \nu_k^{(n)} (f) 
	&= &
	( \omega_k^{(n)})^{-1} \int_0^1 (B_n f(t) -L f(t) ) J_{k-2} (t) dt , \,
	  k=2, \dots , n  , \, f \in C[0,1] .
\end{eqnarray*}
Consequently, for all $f \in C[0,1]$, $x \in [0,1]$,
\begin{eqnarray*}
	 F_n f(x)
	& = &
	 L f(x) 
	+ \sum_{k=2}^n x(1-x) J_{k-2} (x ) \frac{(n-k)! (n+k-1)!}{(n-1)! n!}
	\\
	& &	
	\qquad \times  \int_0^1  (B_n^2 f(t) -L f(t) ) J_{k-2} (t) dt .
\end{eqnarray*}
For each fixed $k$, the corresponding summand tends uniformly to
$$
	x(1-x)  J_{k-2} (x )  \int_0^1 ( f(t) -L f(t) ) J_{k-2} (t) dt , \mbox { when } n \to \infty .
$$ 
It remains to investigate the behavior of the sum when $n \to \infty$.
\begin{remark}
For each $k \geq 2$ and $f \in C[0,1]$,
$$
	\lim_{n \to \infty}  \nu_k^{(n)} (f) 
	=
	  \int_0^1 ( f(t) -L f(t) ) J_{k-2} (t) dt
$$
and 
$$
	\lim_{n \to \infty}  \nu_0^{(n)} (f) = f(0) , \, \lim_{n \to \infty}  \nu_1^{(n)} (f) = f(1)-f(0) .
$$
This solves the Cooper-Waldron type problem for the functionals $ \nu_k^{(n)}$.
\end{remark}


\newpage
\noindent Heiner Gonska
\\Faculty for Mathematics\\University of Duisburg-Essen\\
Forsthausweg 2\\D-47057 Duisburg, Germany
\\
heiner.gonska@uni-due.de
\\[0.3cm]
Margareta Heilmann
\\ Faculty of Mathematics and Natural Sciences\\
University of Wuppertal\\
Gau{\ss}stra{\ss}e 20\\
D-42119 Wuppertal, Germany
\\heilmann@math.uni-wuppertal.de
\\[0.3cm]
Ioan Ra\c{s}a\\
Department of Mathematics\\
Technical University\\
Str. C. Daicoviciu, 15\\
RO-400020 Cluj-Napoca,
Romania\\
Ioan.Rasa@math.utcluj.ro

\section*{Appendix}
\begin{figure}[h]
\includegraphics[scale=0.5]{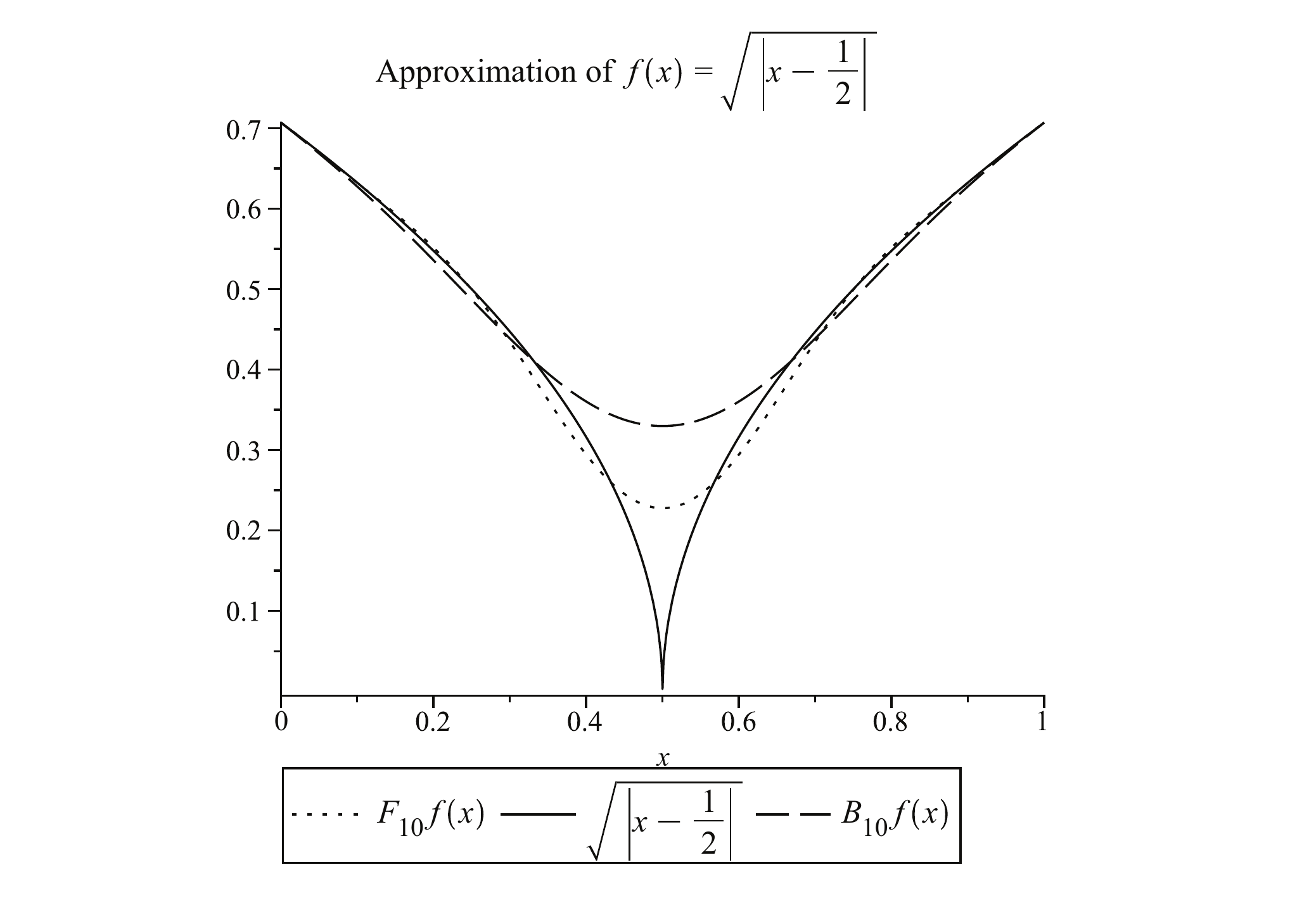}
\end{figure}
\begin{figure}[h]
\includegraphics[scale=0.5]{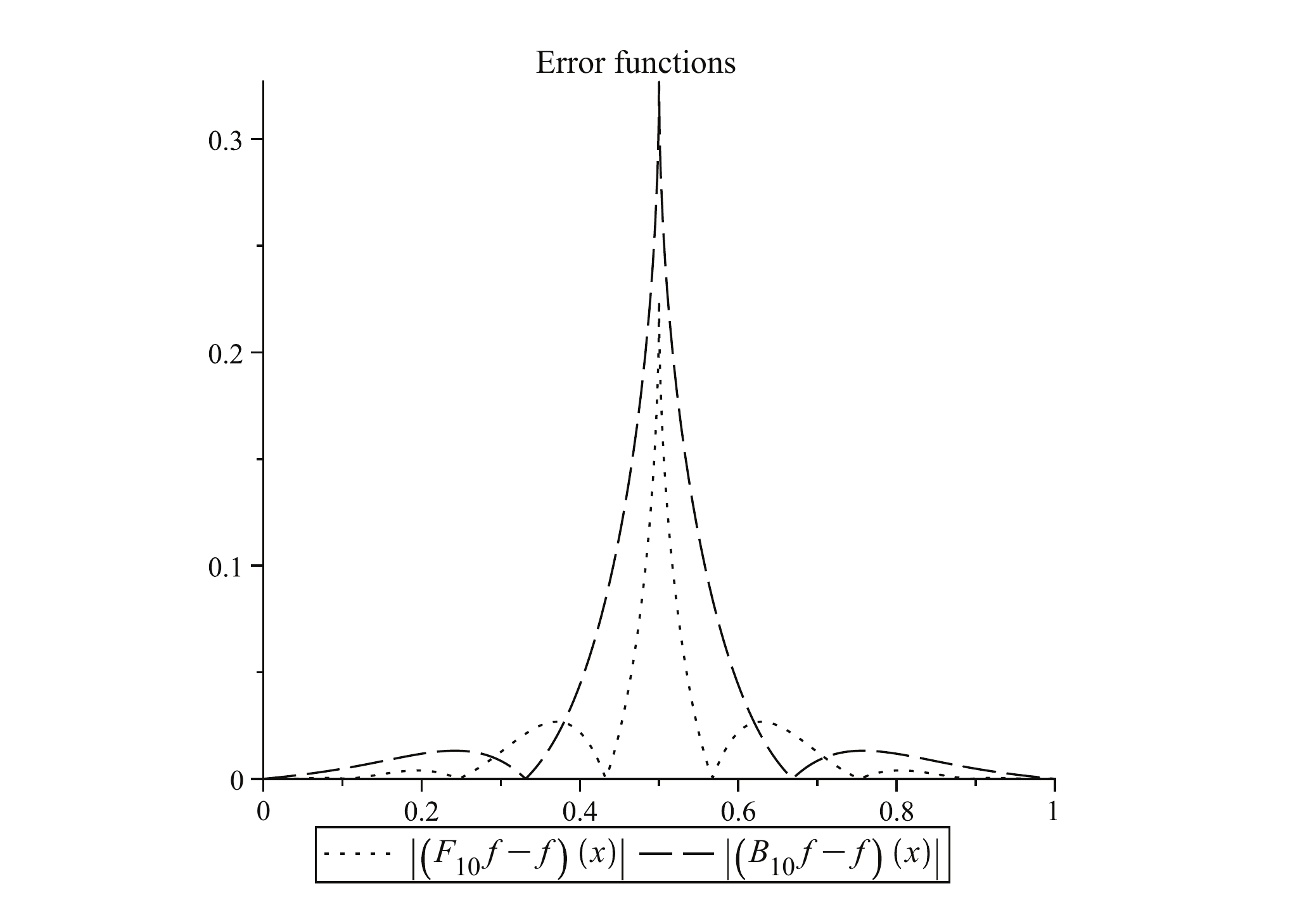}
\end{figure}
\begin{figure}[h]
\includegraphics[scale=0.5]{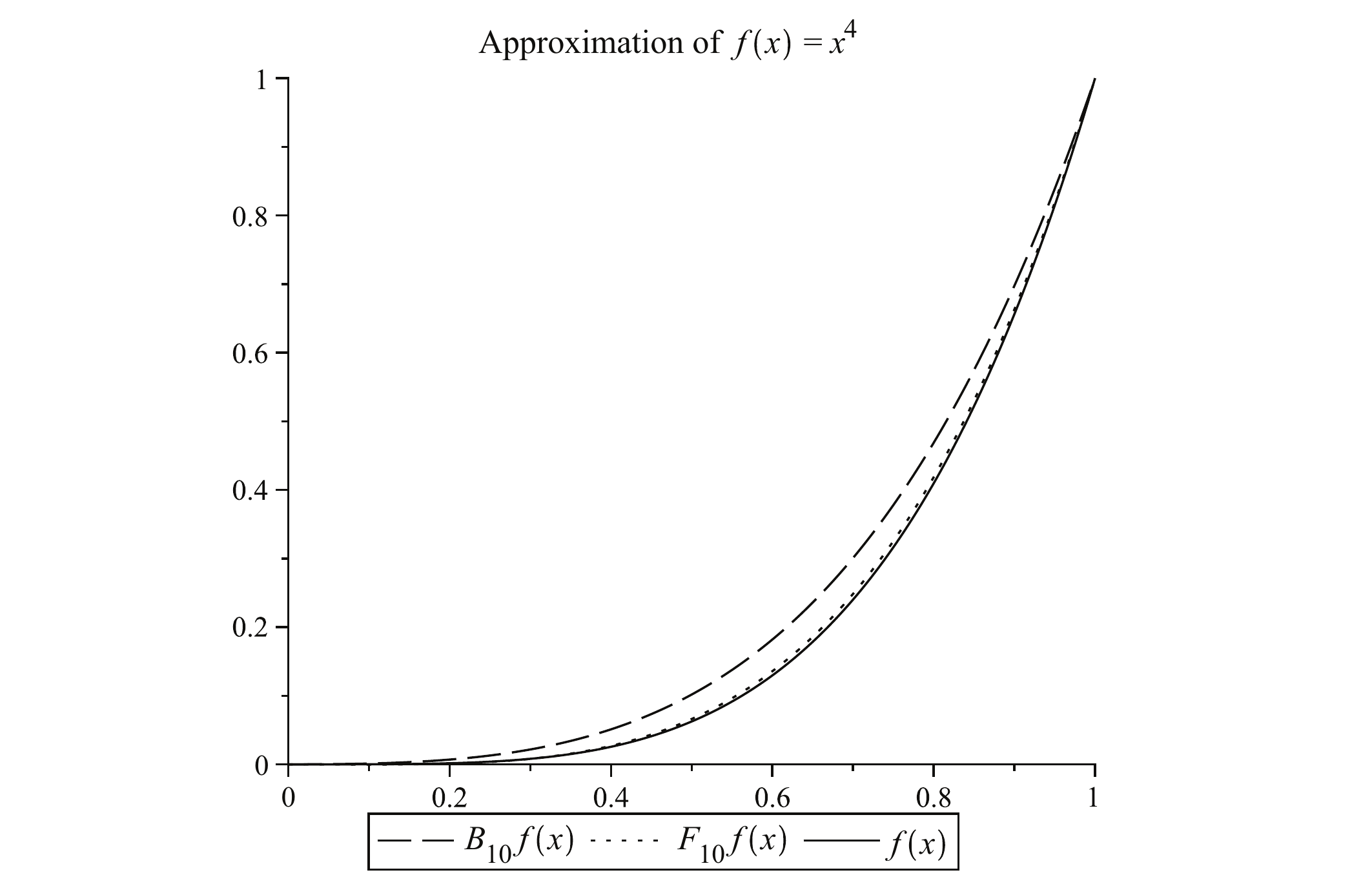}
\end{figure}
\begin{figure}[h]
\includegraphics[scale=0.5]{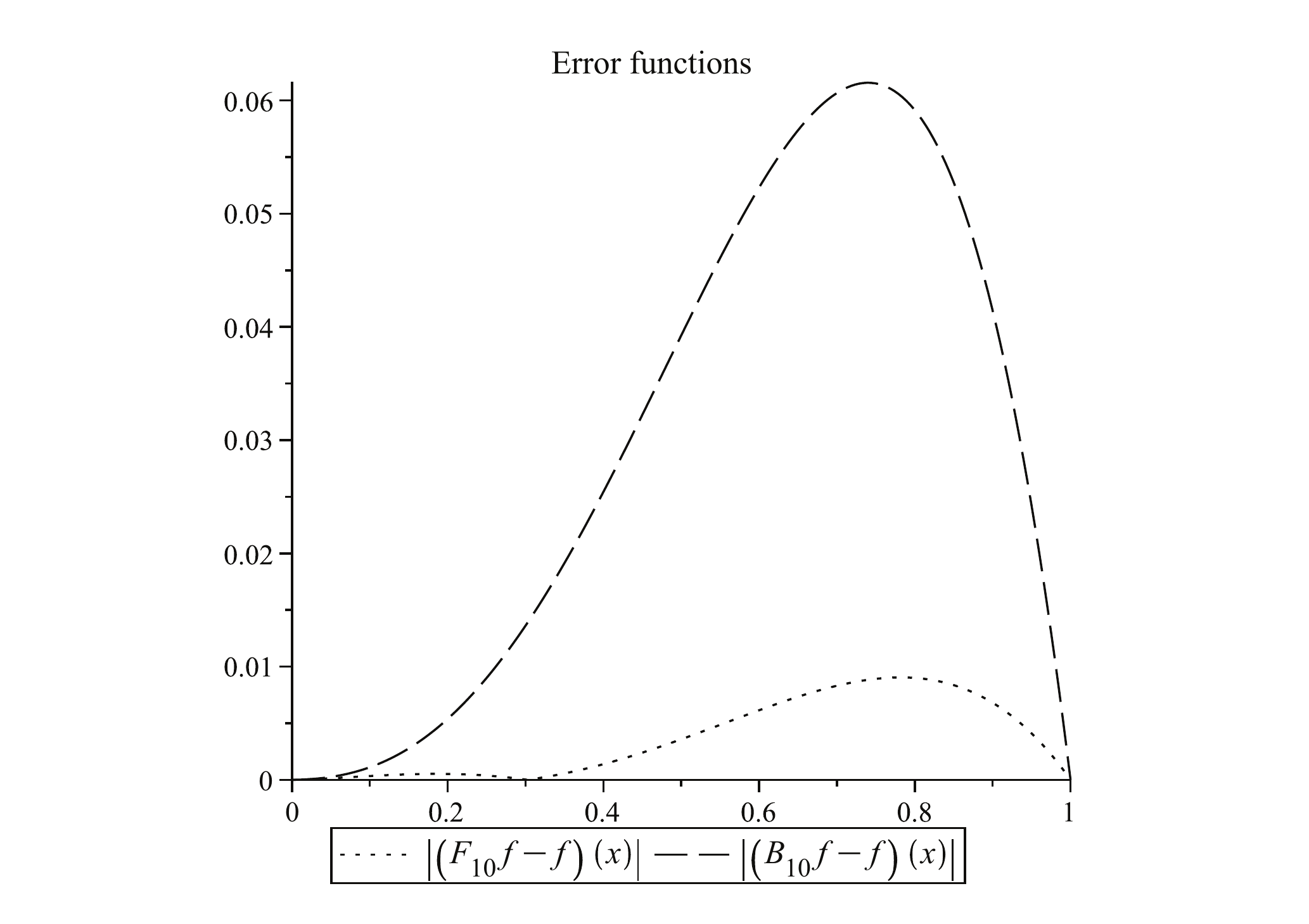}
\end{figure}
\begin{figure}[h]
\includegraphics[scale=0.5]{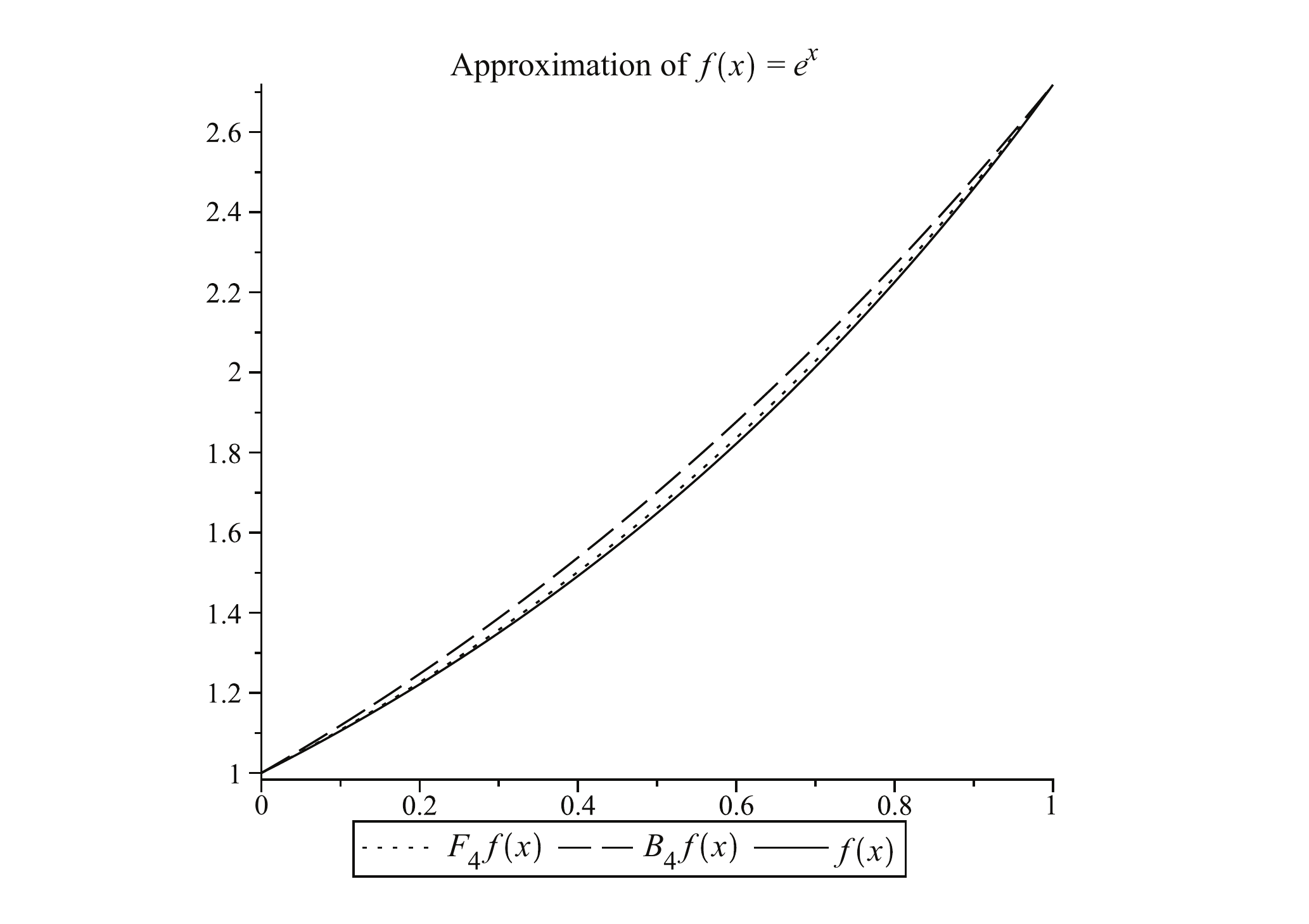}
\end{figure}
\begin{figure}[h]
\includegraphics[scale=0.5]{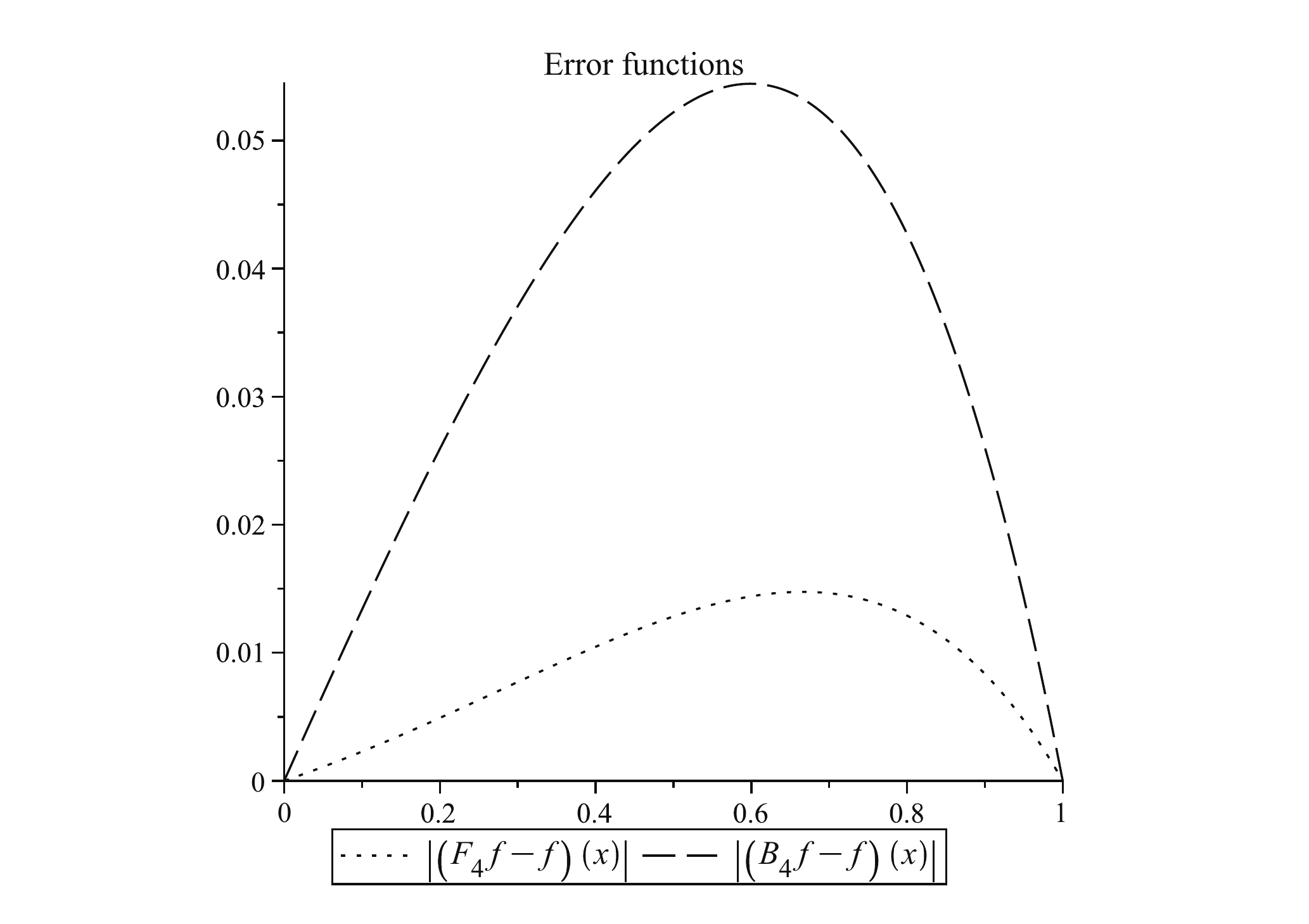}
\end{figure}
\begin{figure}[h]
\includegraphics[scale=0.5]{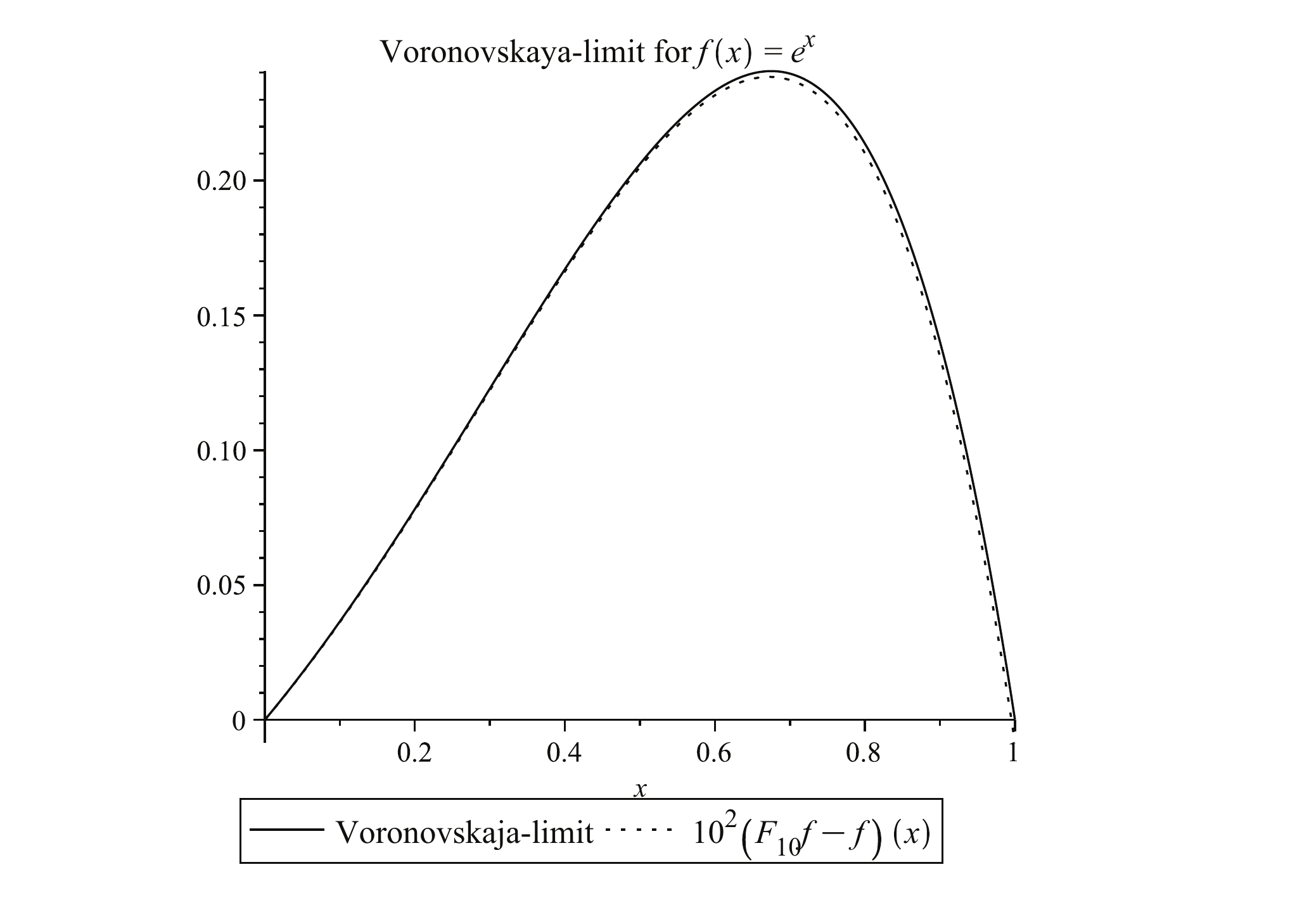}
\end{figure}
\begin{figure}[h]
\includegraphics[scale=0.38]{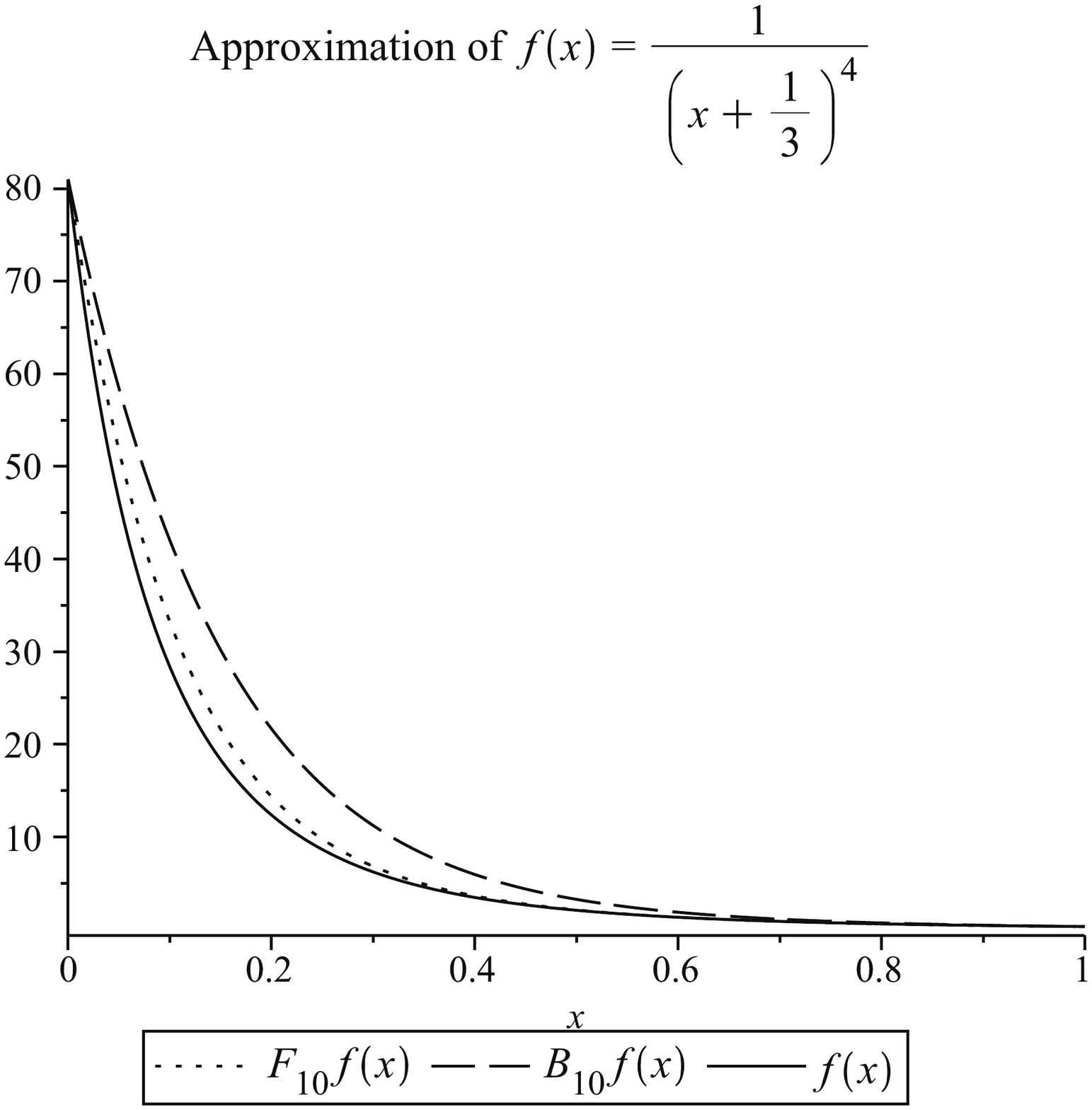}
\end{figure}
\begin{figure}[h]
\includegraphics[scale=0.5]{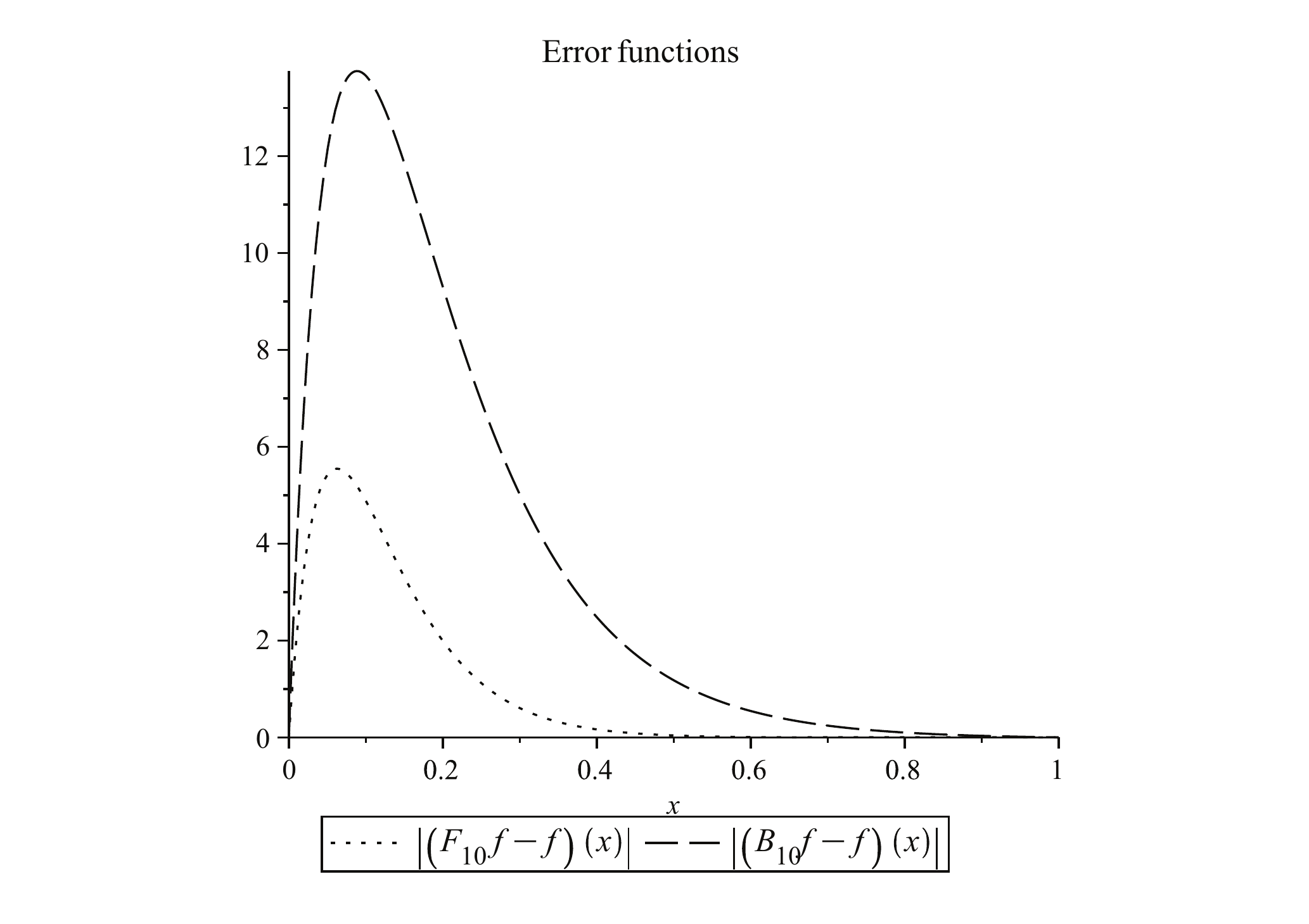}
\end{figure}
\begin{figure}[h]
\includegraphics[scale=0.5]{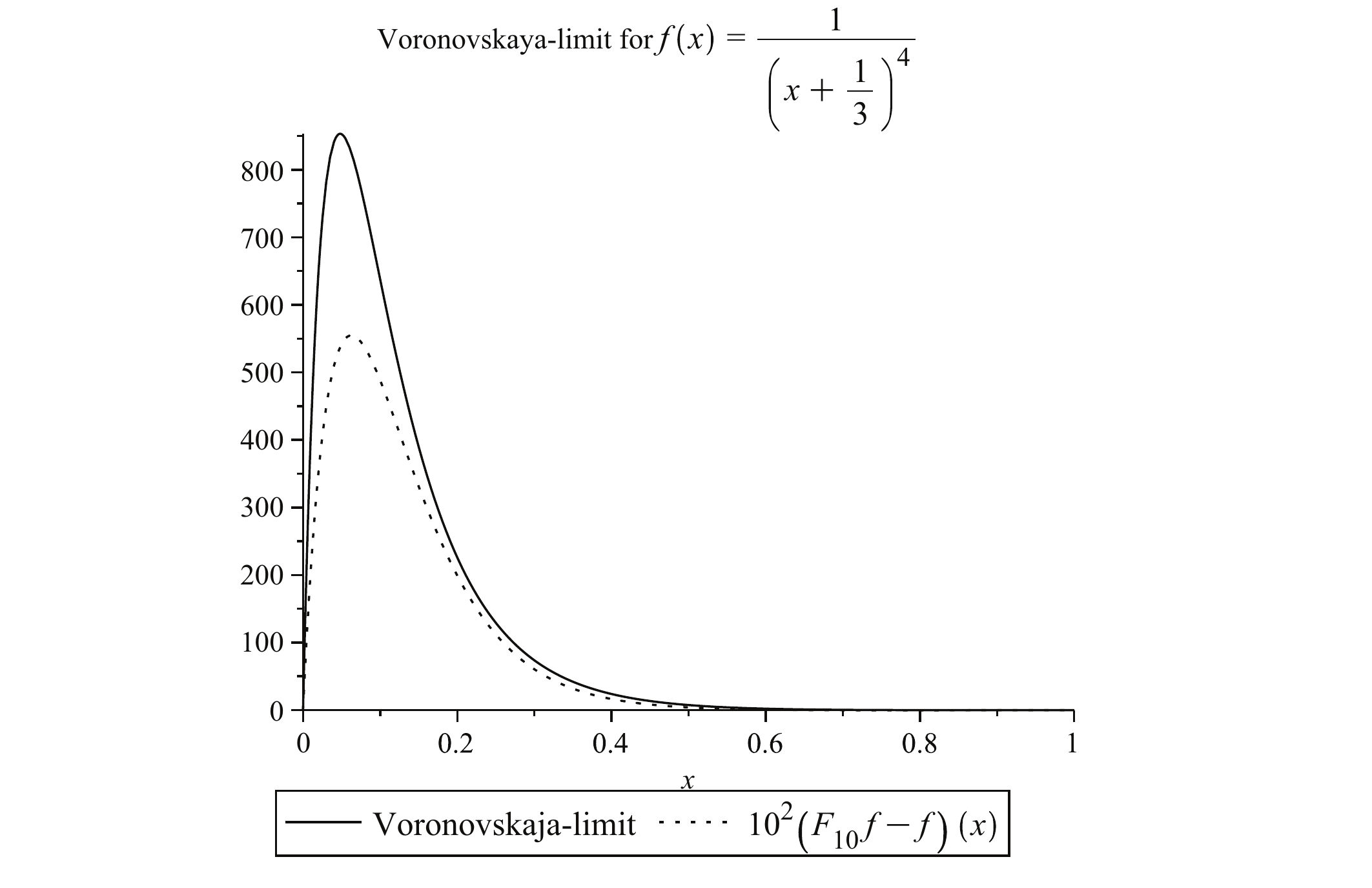}
\end{figure}

\end{document}